\documentclass{article}
\usepackage[utf8]{inputenc}
\usepackage{microtype}
\usepackage{graphicx}
\usepackage{caption}
\usepackage{subcaption}
\usepackage{booktabs}
\usepackage{hyperref}
\usepackage{xcolor}
\usepackage{gensymb}
\usepackage{textcomp}
\usepackage{geometry}

\def\namedlabel#1#2{\begingroup
	#2%
	\def\@currentlabel{#2}%
	\phantomsection\label{#1}\endgroup
}

\newcommand{\footremember}[2]{%
	\footnote{#2}
	\newcounter{#1}
	\setcounter{#1}{\value{footnote}}%
}
\newcommand{\footrecall}[1]{%
	\footnotemark[\value{#1}]%
}

\usepackage{amsmath}
\usepackage{multirow}
\usepackage{amsfonts}
\usepackage{float}
\usepackage{enumitem, hyperref}
\newcommand{\diag}{{\rm diag}}

\def\N{\mathbb{N}}
\newtheorem{eg}{Example}

\newtheorem{prop}{Proposition}
\newtheorem{df}{Definition}
\newtheorem{remark}{Remark}

\author{
	CHEN, Tong \footremember{laas}{LAAS-CNRS, BP 54200, 7 avenue du Colonel Roche, 31031 Toulouse, C\'edex 4, France.}  \\
	\texttt{tchen@laas.fr}
	\and
	LASSERRE, Jean-Bernard \footrecall{laas} \footremember{imt}{IMT, Universit\'e Toulouse 3 Paul Sabatier.}\\
	\texttt{lasserre@laas.fr}
	\and
	MAGRON, Victor \footrecall{laas} \footrecall{imt}\\
	\texttt{vmagron@laas.fr}
	\and
	PAUWELS, Edouard \footremember{irit}{IRIT, Universit\'e de Toulouse, CNRS.} \footrecall{imt} \\
	\texttt{edouard.pauwels@irit.fr}
}

\title{{A Sublevel Moment-SOS Hierarchy for Polynomial Optimization}}

\begin{document}
	\maketitle
	
	\begin{abstract}
		We introduce a sublevel Moment-SOS hierarchy where 	each SDP relaxation can be viewed as an intermediate (or interpolation) between	the $d$-th and $(d+1)$-th order SDP relaxations of the		Moment-SOS hierarchy (dense or sparse version).
		
		With the flexible choice of determining the size (level) and number (depth) of subsets in the SDP relaxation, one is able to obtain different improvements compared to the $d$-th order relaxation, based on the machine memory capacity. In particular, we provide numerical experiments for $d=1$
		and various types of problems both in combinatorial optimization (Max-Cut, Mixed Integer Programming) and deep learning (robustness certification, Lipschitz constant of neural networks), where the standard Lasserre's relaxation (or its sparse variant) is computationally intractable. In our
		numerical results, the lower bounds from the sublevel relaxations improve the bound from Shor's relaxation (first order Lasserre's relaxation) and are significantly closer to the optimal value or to the best-known lower/upper bounds.
	\end{abstract}
	
	\section{Introduction}
	Consider the \textit{polynomial optimization problem (POP)} of the following form:
	\begin{align} \label{pop}
		& f^* := \inf_{\mathbf{x} \in \mathbb{R}^n} \{f(\mathbf{x}): g_i(\mathbf{x}) \ge 0, i = 1, \ldots, p\}\,, \tag{POP}
	\end{align}
	where $f$ and $g_i$ are polynomials in variable $\mathbf{x}$ for all $i = 1, \ldots, p$. 
	Lasserre's hierarchy \cite{lasserre2001global} is a well-known method based on \textit{semidefinite programming (SDP)} to approximate the optimal value of \eqref{pop}, by solving a sequence of SDPs that provide a series of lower bounds and converges to the optimal value of the original problem. Under certain assumptions, such convergence is shown to be finite \cite{nie2014optimality}. 
	
	\paragraph{Related works} 
	{Other related frameworks of relaxations, including  DSOS \cite{ahmadi} based on \textit{linear programming (LP)}, SDSOS \cite{ahmadi} based on \textit{second-order cone programming (SOCP)}, and the hybrid BSOS \cite{lasserre2017bounded} combining the features of LP and SDP hierarchies, also 
		provide lower bounds converging to the optimal value of a POP. 
		Generally speaking, when comparing LP and SDP solvers, the former can handle problems of much larger size. On the other hand, the bounds from LP relaxations 
		are significantly weaker than those obtained by SDP relaxations}, in particular for combinatorial problems \cite{laurent2003comparison}. 
	Based on the standard Lasserre's hierarchy, several further works have
		explored various types  of sparsity patterns inside POPs to compute lower bounds more efficiently and handle larger-scale POPs. The first such extension 
		can be traced back to Waki \cite{waki2006sums} and Lasserre \cite{lasserre2006convergent}
		where the authors consider the so-called \textit{correlative sparsity pattern} (CSP) with
		associated CSP graph whose nodes consist of the POP's variables. Two nodes in the CSP graph are connected via an edge if the two corresponding variables appear in the same constraint or in same monomial of the objective.
		The standard sparse Lasserre's hierarchy splits the full moment and localizing matrices into several smaller submatrices, according to subsets of nodes (\emph{maximal cliques}) in a chordal extension of the CSP graph associated with the POP. When the size of the largest  clique (a crucial parameter of the sparsity pattern) is reasonable 
			the resulting SDP relaxations become tractable.
		There are many successful applications of the resulting \emph{sparse moment-SOS hierarchy}, including certified roundoff error bounds \cite{magron2017certified, magron2018interval},  optimal power flow \cite{josz2018lasserre}, volume computation of sparse semialgebraic sets \cite{tacchi2019exploiting}, approximating regions of attractions of sparse polynomial systems \cite{schlosser2020sparse,tacchi2019approximating}, noncommutative POPs \cite{klep2019sparse}, sparse positive definite functions \cite{mai2020sparse}. Similarly, the sparse BSOS hierarchy \cite{weisser2018sparse} is a sparse version of BSOS for large scale polynomial optimization. 
	
	Besides correlative sparsity, recent developments   \cite{wang2019tssos,wang2020chordal} exploit the so-called \textit{term sparsity}   (TSSOS)  or combine  correlative sparsity and term sparsity (CS-TSSOS) \cite{wang2020cs} to handle large scale polynomial optimization problems. 
	The TSSOS framework relies on a \emph{term sparsity pattern (TSP) graph} whose nodes consist of monomials of some monomial basis.
	Two nodes in a TSP graph are connected via an edge when the product of the corresponding monomials appears in the supports of polynomials involved in the POP or is a monomial of even degree. 
	Extensions have been provided to compute more efficiently approximations of joint spectral radii  \cite{wang2020sparsejsr} and minimal traces  or eigenvalue of noncommutative polynomials \cite{wang2020exploiting}.
	More variants of the sparse moment-SOS hierarchy have been built for quantum bounds on Bell inequalities \cite{P_l_2009}, condensed-matter ground-state problems \cite{Barthel_2012}, quantum many-body problems \cite{haim2020variationalcorrelations}, where one selects a certain subset of words (noncommutative monomials) to decrease the number of SDP variables.
	
	Recently, in \cite{campos2020partial} the authors proposed a \emph{partial} and \emph{augmented partial relaxation} tailored to the Max-Cut problem. It strengthens Shor relaxation by adding some (and not all) constraints from the second-order Lasserre's hierarchy. The same idea was already used in the \emph{multi-order} SDP relaxation of \cite{josz2018lasserre} for solving large-scale optimal power flow (OPF) problems.
		The authors set a threshold for the maximal cliques and include
		the second-order relaxation constraints for the cliques with size under the threshold and 
		the first-order relaxation constraints for the cliques with size over the threshold.
	
	\paragraph{Contribution}
		This work is in the line of research concerned with extensions and/or variants of 
		the Moment-SOS hierarchy so as to handle large-scale POPs out of reach by the standard
		hierarchy.  We provide a principled way to obtain 
		intermediate alternative SDP relaxations between the first- and second-order SDP relaxations of 
		the Moment-SOS hierarchy for general POPs. It encompasses the above cited works 
		\cite{josz2018lasserre,campos2020partial} as special cases for MAX-Cut and OPF problems.
		It can also be generalized to provide intermediate alternative SDP relaxations 
		between (arbitrary) order-$d$ and order-$d+1$ relaxations of the 
		Moment-SOS hierarchy when the order-$d+1$ relaxation is too costly to implement.
	
	
	We develop what we call the \emph{sublevel hierarchy} based on the standard 
		Moment-SOS hierarchy. 
	Compared with existing sparse variants of the latter, 
	we propose several possible SDP relaxations to improve lower bounds for general POPs.

	\paragraph{The basic principle} is quite simple. In the sublevel hierarchy concerned with $d$-th and $(d+1)$-th orders
		of the sparse Moment-SOS hierarchy, from the maximal cliques of a chordal extension of the csp graph, we further select several subsets of nodes (variables). Then in the 
		$d$-th sparse SDP relaxation we also include $(d+1)$-th order moment and localizing matrices w.r.t. these subsets only.
	This methodology reveals helpful if the bound obtained by the $d$-th order relaxation of a POP is not satisfactory  and if one is not able to solve the $(d+1)$-th order relaxation.
	
	One important  distinguished feature of the sublevel hierarchy is to \emph{not} be restricted to POPs with a correlative sparsity pattern. Indeed it can also be applied  to dense POPs or nearly-dense POPs where the problem is sparse except that there are a few dense constraints. As a result we are able to improve bounds obtained at the first-order relaxation (also called Shor's relaxation).
		In  \cite{chen2020semialgebraic} we proposed a heuristic method to deal with general nearly-dense POPs as a trade-off between the first-order and second-order relaxations of the Moment-SOS hierarchy. As we will see
		the heuristic \cite{chen2020semialgebraic} is also a special case of the sublevel hierarchy.
	
	Another feature of the sublevel hierarchy is that we have more flexible ways to tune the resulting relaxation instead of simply increasing the relaxation order (a rigid and costly strategy). More specifically, there are two hyper-parameters in the sublevel hierarchy: (i) the \emph{size} of the selected subsets, and (ii) the \emph{number} of such subsets. 
		Suppose $m$ is the size of a maximal clique of a chordal extension of the csp graph.
	Then we can choose $q$ (called the depth) many subsets of size $l$ (called the level) with $1 \le l \le m$ and $q \le {m \choose l}$. 
	For each maximal clique, we have a wide range of choices for the level and depth, yielding a good trade-off between the solution accuracy and the computational efficiency. 
	
	The outline of the paper is as follows: Section \ref{dense} introduces some preliminaries of dense and sparse Lasserre's hierarchy; Section \ref{sub} is the theoretical part of the sublevel hierarchy and the sublevel relaxation; Section \ref{app} explicitly illustrates the sublevel relaxation for several type of optimization problems; Section \ref{num} shows the results of sublevel relaxation applied on the problems discussed in Section \ref{app}.

	%
	
	\section{Preliminary background on Lasserre's hierarchy} \label{dense}
	In this section we briefly introduce the Lasserre's hierarchy \cite{lasserre2001global} which has already many successful applications in and outside optimization \cite{lasserre2015introduction}.
	First let us recall some notations in polynomial optimization. Given a positive integer $n \in \mathbb{N}$, let $\mathbf{x} = [x_1, \ldots, x_n]^T$ be a vector of decision variables and $\mathbb{R} [\mathbf{x}]$ be the space of real polynomials in variable $\mathbf{x}$. For a set $I \subseteq \{1, 2, \ldots, n\}$, let $\mathbf{x}_{I} := [x_i]_{i \in {I}}$ and let $\mathbb{R} [\mathbf{x}_I]$ be the space of real polynomials in variable $\mathbf{x}_I$. Denote by $\mathbb{R} [\mathbf{x}]$ (resp. $\mathbb{R} [\mathbf{x}]_d$) the vector space of polynomials (resp. of degree at most $d$) in variable $\mathbf{x}$; $\mathcal{P} [\mathbf{x}] \subseteq \mathbb{R} [\mathbf{x}]$ (resp. $\mathcal{P}_d [\mathbf{x}] \subseteq \mathbb{R} [\mathbf{x}]_{2d}$) the convex cone of nonnegative polynomials (resp. nonnegative polynomials of degree at most $2d$) in variable $\mathbf{x}$; $\Sigma [\mathbf{x}] \subseteq \mathcal{P} [\mathbf{x}]$ (resp. $\Sigma [\mathbf{x}]_d \subseteq \mathcal{P}_d [\mathbf{x}]$) the convex cone of SOS polynomials (resp. SOS polynomials of degree at most $2d$) in variable $\mathbf{x}$.
	
	In the context of optimization, Lasserre's hierarchy  allows one to approximate the \emph{global} optimum of \eqref{pop}, by solving a hierarchy 
	of SDPs of increasing size. 
	Each SDP is a semidefinite relaxation of \eqref{pop} in the form:
	\begin{align} \label{mom_opt}
		& \rho_d^{\text{dense}}=\inf_{\mathbf{y}} \{\,L_{\mathbf{y}} (f): L_{\mathbf{y}}(1) = 1, \mathbf{M}_d (\mathbf{y}) \succeq 0, \nonumber\\
		& \qquad\qquad\qquad\quad\mathbf{M}_{d-\omega_i} (g_i \mathbf{y}) \succeq 0, i = 1, \ldots, p\} \,, \tag{Mom-$d$}
	\end{align}
	where $\omega_i = \lceil \deg (g_j)/2 \rceil$, $\mathbf{y}=(y_\alpha)_{\alpha\in\N^n_{2d}}$, $L_{\mathbf{y}}:\mathbb{R}[\mathbf{x}]\to\mathbb{R}$ is the so-called \emph{Riesz linear functional}:
	\[f \:=\sum_\alpha f_\alpha\,\mathbf{x}^\alpha \mapsto L_{\mathbf{y}}(f)\,:=\,\sum_\alpha\,f_\alpha\,y_\alpha,\quad f\in\mathbb{R}[\mathbf{x}],\]
	and $\mathbf{M}_d (\mathbf{y})$, $\mathbf{M}_{d-\omega_i} (g_i \mathbf{y})$ are \emph{moment matrix} and \emph{localizing matrix} respectively; see \cite{lasserre2015introduction} for precise definitions and more details. 
	The semidefinite program \eqref{mom_opt} is the $d$-th order \emph{moment relaxation} of problem \eqref{pop}. As a result, when the semialgebraic set $\mathbf{K}:=\{\mathbf{x}:g_i(\mathbf{x})\geq0,\:i = 1, \ldots, p\}$ is compact, one obtains a monotone sequence of lower bounds $(\rho_d)_{d\in\N}$ with the property $\rho_d\uparrow f^*$ as $d\to\infty$
	under a certain technical \emph{Archimedean} condition; the latter is easily satisfied by including a redundant quadratic  constraint $M-\Vert\mathbf{x}\Vert^2\geq0$ for some well-chosen $M > 0$ in the definition of $\mathbf{K}$ (redundant as $\mathbf{K}$ is compact and $M$ is large enough). At last but not least and interestingly, generically the latter convergence is \emph{finite} \cite{nie2014optimality}. Ideally, one expects an optimal solution $\mathbf{y}^*$ of \eqref{mom_opt} to be the vector of moments up to order $2d$ of the Dirac measure $\delta_{\mathbf{x}^*}$ at a global minimizer $\mathbf{x}^*$ of \eqref{pop}. 

	The hierarchy \eqref{mom_opt} is often referred to as \emph{dense} Lasserre's hierarchy since we do not exploit any possible sparsity pattern of the POP. Therefore, if one solves \eqref{mom_opt} with interior point methods (as  current SDP solvers usually do), then the dense hierarchy 
	is limited to POPs of modest size. 
	Indeed the $d$-th order dense moment relaxation \eqref{mom_opt} involves ${n+2d\choose 2d}$ variables and a moment matrix $M_d(\mathbf{y})$ of size ${n+d\choose d} = O( n^d)$ at fixed $d$. Fortunately, large-scale POPs often exhibit some structured sparsity patterns which can be exploited to yield a \emph{sparse} version of \eqref{mom_opt}, as initially demonstrated in \cite{waki2006sums}. As a result, wider applications of Lasserre's hierarchy have been possible.

	Assume that the set of variables in \eqref{pop} can be divided into $r$ several subsets indexed by $I_k$, for $k \in \{1, \ldots, r\}$, i.e., $\{1, \ldots, n\} = \cup_{k = 1}^r I_k$. Suppose that the following assumptions hold:
	
	\textbf{A1}: The function $f$ is a sum of polynomials,  each 
	summand	involving variables of only one subset, i.e., $f(\mathbf{x}) = \sum_{k = 1}^r f_k(\mathbf{x}_{I_{k}})$;
	
	\textbf{A2}: Each constraint also involves variables of only one subset, i.e., $g_i \in \mathbb{R} [\mathbf{x}_{I_{k(i)}}]$ for some $k(i) \in \{1, \cdots, r\}$;
	
	\textbf{A3}: The subsets $I_k$ satisfy the \emph{Running Intersection Property (RIP)}: for every $k \in \{1, \cdots, r-1\}$, $I_{k+1} \cap \bigcup_{j = 1}^k I_j \subseteq I_s$, for some $s \le k$.
	
	It turns out that the maximal cliques in the chordal extension of the csp graph induced by the POP satisfy the RIP \cite{waki2006sums}. From now on, we will call the these subsets \emph{cliques}, in order to distinguish from the subsets in the sublevel hierarchy that will be discussed in the next section. A POP with a sparsity pattern is of the form:
	\begin{align} \label{sparse_opt}
		\inf_{\mathbf{x} \in \mathbb{R}^n} \{f(\mathbf{x}): g_i(\mathbf{x}_{I_{k}}) \ge 0, 
		\:i = 1, \ldots, p; \: i \in I_k\}\,, \tag{SpPOP}
	\end{align}
	and its associated \emph{sparse Lasserre's hierarchy} reads:
	\begin{align} \label{sparse_mom_opt}
		& \rho_d^{\text{sparse}}=\inf_{\mathbf{y}} \{L_{\mathbf{y}} (f): L_{\mathbf{y}}(1) = 1, \mathbf{M}_d (\mathbf{y}, I_k) \succeq 0, k \in \{1, \cdots, r\}, \nonumber\\
		& \qquad\qquad\qquad\quad\mathbf{M}_{d-\omega_i} (g_i \,\mathbf{y}, I_{k}) \succeq 0\,,\: i \in \{1, \cdots, p\}; \: i \in I_k\,\}\,, \tag{SpMom-$d$}
	\end{align}
	where $d$, $\omega_i$, $\mathbf{y}$, $L_{\mathbf{y}}$ are defined as in \eqref{mom_opt}
	but with a crucial difference. The matrix $\mathbf{M}_d(\mathbf{y},I_k)$ (resp. $\mathbf{M}_{d-\omega_i}(g_i\,\mathbf{y},I_k)$) is a submatrix of the moment matrix $\mathbf{M}_d (\mathbf{y})$ (resp. localizing matrix $\mathbf{M}_{d-\omega_i} (g_i \mathbf{y})$) with respect to the clique $I_k$, and hence of much smaller size
	${\tau_k+d\choose \tau_k}$ if $\vert I_k\vert=:\tau_k\ll n$. Finally, $\rho_d^{\text{sparse}} \leq f^*$ for all $d$ and moreover, if the cliques $I_k$ satisfy the RIP, then we still obtain the convergence $\rho_d^{\text{sparse}} \uparrow f^*$ as $d\to\infty$, as for the dense relaxation \eqref{mom_opt}.
	
	Finally, for each fixed $d$, the dual of \eqref{mom_opt}
	reads:
	\begin{align*} \label{sos_opt}
		\sup_{t \in \mathbb{R}} \{t: f - t = \theta + \sum_{i = 1}^p \sigma_i g_i\} \,, \tag{SOS-$d$}
	\end{align*}
	where $\theta$ is a \emph{sum-of-squares (SOS)} polynomial in $\mathbb{R} [\mathbf{x}]$ of degree at most $2d$, and $\sigma_j$ are SOS polynomials in $\mathbb{R} [\mathbf{x}]$ of degree at most $2(d-\omega_i)$ with $\omega_i = \lceil \deg (g_j)/2 \rceil$. The right-hand-side of the identity in \eqref{sos_opt} is nothing less than 
	Putinar's positivity certificate \cite{putinar1993positive} for the polynomial 
	$\mathbf{x}\mapsto f(\mathbf{x})-t$	on the compact semialgebraic set $\mathbf{K}$.

	Similarly, the dual problem of \eqref{sparse_mom_opt} reads:
	\begin{align*} \label{sparse_sos_opt}
		& \sup_{t \in \mathbb{R}} \{t: f - t = \sum_{k = 1}^m \big(\theta_{k} + \sum_{i \in I_k} \sigma_{i,k} g_i\big)\} \,, \tag{SpSOS-$d$}
	\end{align*}
	where $\theta_{k}$ is an SOS in $\mathbb{R} [\mathbf{x}_{I_k}]$ of degree at most $2d$, and $\sigma_{i,k}$ is an SOS in $\mathbb{R} [\mathbf{x}_{I_k}]$ of degree at most $2(d-\omega_i)$ with $\omega_i = \lceil \deg (g_i)/2 \rceil$, for each $k=1,\dots,p$. Then \eqref{sparse_sos_opt}
	implements the sparse Putinar's positivity certificate \cite{lasserre2006convergent, waki2006sums}.
	
	\begin{eg} \label{eg}
		Let $\mathbf{x} \in \mathbb{R}^6, \mathbf{x}_{1:4} := [x_i]_{i = 1}^4, \mathbf{x}_{3:6} := [x_i]_{i = 3}^6$. 
		We minimize $f (\mathbf{x}) = - ||\mathbf{x}||_2^2$, under the semialgebraic set defined by $g_1 (\mathbf{x}) = 1 - ||\mathbf{x}_{1:4}||_2^2 \ge 0$ and $g_2 (\mathbf{x}) = 1 - ||\mathbf{x}_{3:6}||_2^2 \ge 0$. Then, the second-order dense Lasserre's relaxation reads
		$$\sup_{t \in \mathbb{R}} \{t: f(\mathbf{x}) - t = \theta (\mathbf{x}) + \sigma_1 (\mathbf{x}) g_1 (\mathbf{x}) + \sigma_2 (\mathbf{x}) g_2 (\mathbf{x})\}$$
		where $\theta$ is a degree-4 SOS polynomial in variable $\mathbf{x}$, $\sigma_1, \sigma_2$ are degree-2 SOS polynomials in variable $\mathbf{x}$. Define $I_1 = \{1,2,3,4\}$ and $I_2 = \{3,4,5,6\}$, then $g_1 \in \mathbb{R} [\mathbf{x}_{I_1}]$ and $g_2 \in \mathbb{R} [\mathbf{x}_{I_2}]$. The second-order sparse Lasserre's relaxation reads
		$$\sup_{t \in \mathbb{R}} \{t: f(\mathbf{x}) - t = \big(\theta_{1} (\mathbf{x}_{I_1}) + \sigma_1 (\mathbf{x}_{I_1}) g_1 (\mathbf{x})\big) + \big(\theta_{2} (\mathbf{x}_{I_2}) + \sigma_2 (\mathbf{x}_{I_2}) g_2 (\mathbf{x})\big)\}$$
		where $\theta_{k}$ is a degree-4 SOS polynomials in variable $\mathbf{x}_{I_k}$, and  $\sigma_k$ is a degree-2 SOS polynomials in variable $\mathbf{x}_{I_k}$, for eac $k=1,2$.
	\end{eg}
	
	\section{Sublevel hierarchy}\label{sub}
	As seen in Section \ref{dense}, the way to reduce the size of the moment and localizing matrices in \eqref{mom_opt} is either by reducing the relaxation order or the number of variables/terms in the SOS weights involved in the Putinar's representation. The authors from \cite{josz2018lasserre} propose the \textit{multi-order} Lasserre's hierarchy to deal with large-scale optimal power flow problems. In this hierarchy, one reduces the relaxation order with respect to the constraints with large number of variables. This approach is reused as the so-called \emph{partial} relaxation to solve Max-Cut problems in \cite{campos2020partial}. The authors in \cite{campos2020partial} also proposed the augmented partial relaxation as an extended version of the partial relaxation, to improve the bounds further. In this section, we develop the sublevel hierarchy which is a generalization of several existing frameworks for both sparse and non-sparse POPs, and show that in the case of Max-Cut problems, the partial and augmented partial relaxation can be cast as special instances of the sublevel relaxation.

	\subsection{Deriving the sublevel hierarchy} \label{sub1}
	For problem \eqref{pop}, the $d$-th order dense Lasserre's relaxation relates to the Putinar's certificate $f - t = \sigma_0 + \sum_{i = 1}^p \sigma_i g_i$ where $\sigma_0$ is an SOS in $\mathbb{R} [\mathbf{x}]$ of degree at most $2d$ and $\sigma_i$ are SOS in $\mathbb{R} [\mathbf{x}]$ of degree at most $2(d-\omega_i)$ with $\omega_i = \lceil \deg (g_i) /2 \rceil$. 
	In this section, we are going to choose some subsets of the variable $\mathbf{x}$ to decrease the number of terms involved in the SOS multipliers $\sigma_0$ and $\sigma_i$, and define the intermediate sublevel hierarchies between the $d$-th and $(d+1)$-th order relaxations. 
	
	Note that in the dense variant of Lasserre's hierarchy, one approximates the cone of positive polynomials from the inside with the following hierarchy of SOS cones:
	$$\mathbb{R} = \Sigma [\mathbf{x}]_0 \subseteq \Sigma [\mathbf{x}]_1 \subseteq \ldots \subseteq \Sigma [\mathbf{x}]$$ 
	with $\bigcup_{d = 0}^{+\infty} \Sigma [\mathbf{x}]_d = \Sigma [\mathbf{x}]$. 
	Similarly, in the sparse variant, one relies on the following hierarchy of direct sums of SOS cones:
	$$\mathbb{R} = \oplus_k \Sigma [\mathbf{x}_{I_k}]_0 \subseteq \oplus_k \Sigma [\mathbf{x}_{I_k}]_1 \subseteq \ldots \subseteq \oplus \Sigma [\mathbf{x}_{I_k}]$$
	with $\bigcup_{d = 0}^{+\infty} (\oplus_k \Sigma [\mathbf{x}_{I_k}]_d) = \oplus_k \Sigma [\mathbf{x}_{I_k}]$.
	
	\begin{df}
		\textbf{(Sublevel hierarchy of SOS cones)} Let $n$ be the number of variables in \eqref{pop}. For $d \ge 1$ and $0 \le l \le n$, the $l$-th level SOS cone associated to $\Sigma [\mathbf{x}]_d$, denoted by $\Sigma [\mathbf{x}]_{d}^l$, is an SOS cone lying between $\Sigma[\mathbf{x}]_d$ and $\Sigma[\mathbf{x}]_{d+1}$, which is defined as
		$$\Sigma[\mathbf{x}]_d \subseteq \Sigma[\mathbf{x}]_{d}^l := \Sigma[\mathbf{x}]_d + \tilde{\Sigma}[\mathbf{x}]_{d+1}^l \subseteq \Sigma[\mathbf{x}]_{d+1}$$
		where $\tilde{\Sigma}[\mathbf{x}]_{d+1}^l: = \bigg\{\displaystyle \sum_{|I| = l} \sigma_I (\mathbf{x}_{I}): I \subseteq \{1, \ldots, n\}, \sigma_I (\mathbf{x}_{I}) \in \Sigma [\mathbf{x}_I]_{d+1} \bigg\} \subseteq \Sigma[\mathbf{x}]_{d+1}$, i.e., the SOS polynomials in $\tilde{\Sigma}[\mathbf{x}]_{d+1}^l$ are the elements in $\Sigma[\mathbf{x}]_{d+1}$ which can be decomposed into several components where each component is an SOS polynomial in $l$ variables. Let us use the convention $\Sigma[\mathbf{x}]_{d}^0 := \Sigma[\mathbf{x}]_{d}$. 
		Then, for the dense case, we rely on the sublevel hierarchy of inner approximations of the cone of positive polynomials:
		$$\Sigma [\mathbf{x}]_d = \Sigma [\mathbf{x}]_{d}^0 \subseteq \Sigma [\mathbf{x}]_{d}^1 \subseteq \ldots \subseteq \Sigma [\mathbf{x}]_{d}^n = \Sigma [\mathbf{x}]_{d+1}$$

		Similarly, suppose that $\{I_k\}_{1\leq k\leq r}$ are the cliques of the sparse problem \eqref{sparse_opt}. For $l \le \tau_k := |I_k|$, we define the $l$-th level SOS cone of $\Sigma [\mathbf{x}_{I_k}]_d$, denoted by $\Sigma [\mathbf{x}_{I_k}]_{d}^l$, as
		$$\Sigma[\mathbf{x}_{I_k}]_d \subseteq \Sigma[\mathbf{x}_{I_k}]_{d}^l := \Sigma[\mathbf{x}_{I_k}]_d + \tilde{\Sigma}[\mathbf{x}_{I_k}]_{d+1}^l \subseteq \Sigma[\mathbf{x}_{I_k}]_{d+1}$$
		where $\tilde{\Sigma}[\mathbf{x}_{I_k}]_{d+1}^l: = \bigg\{\displaystyle \sum_{|I| = l} \sigma_I (\mathbf{x}_{I}): I \subseteq I_k, \sigma_I (\mathbf{x}_{I}) \in \Sigma [\mathbf{x}_I]_{d+1} \bigg\} \subseteq \Sigma[\mathbf{x}_{I_k}]_{d+1}$, i.e., the SOS polynomials in $\tilde{\Sigma}[\mathbf{x}_{I_k}]_{d+1}^l$ are the elements in $\Sigma[\mathbf{x}_{I_k}]_{d+1}$ which can be decomposed into several components where each component is an SOS polynomial in $l$ variables indexed by $I_k$. Then, for the sparse case, we rely on the sublevel hierarchy of inner approximations of the cone of positive polynomials:
		$$\Sigma [\mathbf{x}_{I_k}]_d = \Sigma [\mathbf{x}_{I_k}]_{d}^0 \subseteq \Sigma [\mathbf{x}_{I_k}]_{d}^1 \subseteq \ldots \subseteq \Sigma [\mathbf{x}_{I_k}]_{d}^{\tau_k} = \Sigma [\mathbf{x}_{I_k}]_{d+1}$$
	\end{df}
	
	\begin{remark} Lasserre's hierarchy relies on a  hierarchy of SOS cones, while the sublevel hierarchy relies on a hierarchy of sublevel SOS cones. 
		Take the sparse case for illustration, solving the $d$-th order relaxation of the standard sparse Lasserre's hierarchy boils down to finding SOS multipliers in the cone $\Sigma [\mathbf{x}_{I_k}]_d \oplus \Sigma [\mathbf{x}_{I_k}]_{d-\omega_i}$ for each clique $I_k$, i.e., $\bigoplus_k (\Sigma [\mathbf{x}_{I_k}]_d \oplus \Sigma [\mathbf{x}_{I_k}]_{d-\omega_i})$. 
		Solving the $d$-th order sublevel hierarchy boils down to finding SOS multipliers in the intermediate cones $\bigoplus_k (\Sigma [\mathbf{x}_{I_k}]_{d}^{l_k} \oplus \Sigma [\mathbf{x}_{I_k}]_{d-\omega_i}^{l_k})$ for some $0 \le l_k \le \tau_k$. 
		This cone approximates the standard cone $\bigoplus_k (\Sigma [\mathbf{x}_{I_k}]_d \oplus \Sigma [\mathbf{x}_{I_k}]_{d-1})$ as $l_k$ gets larger since $\bigoplus_k (\Sigma [\mathbf{x}_{I_k}]_{d}^{\tau_k} \oplus \Sigma [\mathbf{x}_{I_k}]_{d-\omega_i}^{\tau_k}) = \bigoplus_k (\Sigma [\mathbf{x}_{I_k}]_d \oplus \Sigma [\mathbf{x}_{I_k}]_{d-\omega_i})$. 
		We will see in the next definition that this is the so-called sublevel relaxation, and we call the vector $\{l_k\}$ the vector of sublevels of the relaxation. Each $l_k$ determines the size of the subsets in the clique $I_k$ and is called a \emph{sublevel}.
	\end{remark}
	
	\begin{df}\label{sublevel}
		\textbf{(Sublevel hierarchy of moment-SOS relaxations)} Let $n$ be the number of variables in \eqref{pop}. For each constraint $g_i \ge 0$ in \eqref{pop}, we define a sublevel $0 \le l_i \le n$ and a depth $0 \le q_i \le n$. 
		Denote by $\mathbf{l} = [l_i]_{i = 1}^p$ the vector of sublevels and $\mathbf{q} = [q_i]_{i = 1}^p$ the vector of depths. 
		Then, the $(\mathbf{l}, \mathbf{q})$-sublevel relaxation of the $d$-th order dense SOS problem \eqref{sos_opt} reads
		\begin{align*} \label{sos_opt_sublevel}
			\sup_{t \in \mathbb{R}} \bigg\{t: f - t = \theta_0 + \sum_{i = 1}^p \big(\tilde{\theta}_{i} + (\sigma_i + \tilde{\sigma}_i) g_i\big)\bigg\} \,, \tag{SubSOS-$[d,\mathbf{l}, \mathbf{q}]$}
		\end{align*}
		where $\theta_0$ (resp. $\sigma_i$) are SOS polynomials in $\Sigma [\mathbf{x}]_{d}$ (resp. $\Sigma [\mathbf{x}]_{d - \omega_i}$), and $\tilde{\theta}_{i}$ (resp. $\tilde{\sigma}_i$) are SOS polynomials in $\tilde{\Sigma}[\mathbf{x}]_{d+1}^{l_i}$ (resp. $\tilde{\Sigma}[\mathbf{x}]_{d-\omega_i+1}^{l_i}$) with $\omega_i = \lceil \deg (g_i)/2 \rceil$). Moreover, each $\tilde{\sigma}_i$ is a sum of $q_i$ SOS polynomials where each sum term involves variables in a certain subset $\Gamma_{i,j} \subseteq \{1, 2, \ldots, n\}$ with $|\Gamma_{i,j}| = l_i$, i.e., $\tilde{\sigma}_i = \sum_{j = 1}^{q_i} \tilde{\sigma}_{i,j}$ where $\tilde{\sigma}_{i,j} \in \Sigma [\mathbf{x}_{\Gamma_{i,j}}]_{d-\omega_i+1}$. Each $\tilde{\theta}_{i}$ is also a sum of $q_i$ SOS polynomials where the sum terms share the same variable sets $\Gamma_{i,j}$ as $\tilde{\sigma}_{i,j}$, i.e., $\tilde{\theta}_{i} = \sum_{j = 1}^{q_i} \tilde{\theta}_{i,j}$ where $\tilde{\theta}_{i,j} \in \Sigma [\mathbf{x}_{\Gamma_{i,j}}]_{d+1}$. The equation \eqref{sos_opt_sublevel} can be compressed as an analogical form of the standard dense Lasserre's relaxation:
		\begin{align*}
			\sup_{t \in \mathbb{R}} \bigg\{t: f - t = \sum_{i = 1}^p (\tilde{\theta}_{i} + \tilde{\sigma}_i g_i)\bigg\} \,,
		\end{align*}
		where $\tilde{\theta}_{i}$ (resp. $\tilde{\sigma}_i$) are SOS polynomials in $\Sigma[\mathbf{x}]_{d+1}^l$ (resp. $\Sigma[\mathbf{x}]_{d-\omega_i+1}^l$).
		
		Similarly, suppose that $(I_k)_{1\leq k\leq p}$ are the cliques of the sparse problem \eqref{sparse_opt} with $\tau_k = |I_k|$. For each constraint $g_i \ge 0$ in \eqref{sparse_opt}, denote by $k(i)$ the set of indices $s$ such that $i \in I_{s}$. For each $i$ and $s \in k(i)$, define a sublevel $0 \le l_{i,s} \le \tau_{s}$ and a depth $0 \le q_{i,s} \le \tau_{s}$. Denote by $\mathbf{l} = [l_{i,s}]_{i = 1, \ldots, p; \: s \in k(i)}$ the vector of sublevels and $\mathbf{q} = [q_{i,s}]_{i = 1, \ldots, p; \: s_i \in k(i)}$ the vector of depths. Then, the $(\mathbf{l}, \mathbf{q})$-sublevel relaxation of the $d$-th order sparse SOS problem \eqref{sparse_sos_opt} reads
		\begin{align*} \label{sparse_sos_opt_sublevel}
			& \sup_{t \in \mathbb{R}} \bigg\{t: f - t = \sum_{k = 1}^m \bigg(\theta_{0,k} + \sum_{i \in I_k} \big(\tilde{\theta}_{i,k} + (\sigma_{i,k}+ \tilde{\sigma}_{i,k}) g_i\big)\bigg)\bigg\} \,, \tag{SubSpSOS-$[d, \mathbf{l}, \mathbf{q}]$}
		\end{align*}
		where $\theta_{0,k}$ (resp. $\sigma_{i,k}$) are SOS polynomials in $\Sigma [\mathbf{x}_{I_k}]_{d}$ (resp. $\Sigma [\mathbf{x}_{I_k}]_{d - \omega_i}$), and $\tilde{\theta}_{0,k}$ (resp. $\tilde{\sigma}_{i,k}$) are SOS polynomials in $\tilde{\Sigma}[\mathbf{x}_{I_k}]_{d+1}^{l_{i,k}}$ (resp. $\tilde{\Sigma}[\mathbf{x}_{I_k}]_{d-\omega_i+1}^{l_{i,k}}$) with $\omega_i = \lceil \deg (g_i)/2 \rceil$. Moreover, each $\tilde{\sigma}_{i,k}$ with $i \in I_k$ is a sum of $q_{i,k}$ SOS polynomials where each sum term involves variables in a certain subset $\Gamma_{i,k,j} \subseteq I_k$ with $|\Gamma_{i,k,j}| = l_{i,k}$, i.e., $\tilde{\sigma}_{i,k} = \sum_{j = 1}^{q_{i,k}} \tilde{\sigma}_{i,k,j}$ where $\tilde{\sigma}_{i,k,j} \in \Sigma [\mathbf{x}_{\Gamma_{i,k,j}}]_{d - \omega_i +1}$. Each $\tilde{\theta}_{i,k}$ is also a sum of $q_{i,k}$ SOS polynomials where the sum terms share the same variable sets $\Gamma_{i,k,j}$ as $\tilde{\sigma}_{i,k,j}$, i.e., $\tilde{\theta}_{i,k} = \sum_{j = 1}^{q_{i,k}} \tilde{\theta}_{i,k,j}$ where $\tilde{\theta}_{i,k,j} \in \Sigma [\mathbf{x}_{\Gamma_{i,k,j}}]_{d+1}$. The equation \eqref{sparse_sos_opt_sublevel} can also be compressed as an analogical form of the standard sparse Lasserre's relaxation:
		\begin{align*}
			& \sup_{t \in \mathbb{R}} \bigg\{t: f - t = \sum_{k = 1}^m \sum_{i \in I_k} \big(\tilde{\theta}_{i,k} + \tilde{\sigma}_{i,k} g_i\big)\bigg\} \,,
		\end{align*}
		where $\tilde{\theta}_{i,k}$ (resp. $\tilde{\sigma}_{i,k}$) are SOS polynomials in $\Sigma [\mathbf{x}_{I_k}]_{d+1}^l$ (resp. $\Sigma [\mathbf{x}_{I_k}]_{d - \omega_i +1}^l$).
	\end{df}
	
	\begin{remark}
		(i). If one of the sublevel $l_i$ (resp. $l_{i,k}$) in the dense (resp. sparse) sublevel relaxation is such that $l_i = n$ (resp. $l_{i,k} = \tau_k$), then the depth $q_i$ (resp. $q_{i,k}$) should automatically be 1.
		
		(ii). The heuristics to determine the subsets ($\Gamma_{i,j}$ for the dense case and $\Gamma_{i,k,j}$ for the sparse case) in the sublevel relaxation will be discussed in the next section.
		
		(iii). The size of the SDP Gram matrix associated to an SOS polynomial in $\Sigma [\mathbf{x}]_{d}^l$ (resp. $\Sigma [\mathbf{x}_{I_k}]_{d}^l$) is $\max \{\binom{n+d}{d}, \binom{l+d+1}{d+1} \}$ (resp. $\max \{\binom{|I_k|+d}{d}, \binom{l+d+1}{d+1} \}$). 
		If the lower bound obtained by solving the SOS problem over  $\Sigma [\mathbf{x}]_{d+1}$ (resp. $\Sigma [\mathbf{x}_{I_k}]_{d+1}$) is not satisfactory enough, then we may try to find more accurate solutions in one of the cones of $\Sigma [\mathbf{x}]_{d}^l$ (resp. $\Sigma [\mathbf{x}_{I_k}]_{d}^l$).
	\end{remark}
	
	\begin{eg}
		Take the polynomials $f, g_k$ and the cliques $I_k$ as in Example \ref{eg}. Define $\mathbf{l} = [2,2]$ and $\mathbf{q} = [1,1]$. We select subsets w.r.t. $g_1$ and $g_2$ respectively as $\Gamma_{1,1} = \{1,2\}$, $\Gamma_{2,1} = \{5,6\}$. Then, the second-order dense $(\mathbf{l}, \mathbf{q})$-sublevel relaxation reads
		$$\sup_{t \in \mathbb{R}} \{t: f(\mathbf{x}) - t = \theta_0 (\mathbf{x}) + \big(\tilde{\theta_1} (\mathbf{x}_{\Gamma_{1,1}}) + \tilde{\sigma_1} (\mathbf{x}_{\Gamma_{1,1}}) g_1 (\mathbf{x})\big) + \big(\tilde{\theta_2} (\mathbf{x}_{\Gamma_{2,1}}) + \tilde{\sigma_2} (\mathbf{x}_{\Gamma_{2,1}}) g_2 (\mathbf{x})\big)\}$$
		where $\theta_0$ is a degree-2 SOS polynomial in variable $\mathbf{x}$, $\tilde{\theta_k}$ are degree-4 SOS polynomials in variable $\mathbf{x}_{\Gamma_{k,1}}$, $\tilde{\sigma_k}$ are degree-2 SOS polynomials in variable $\mathbf{x}_{\Gamma_{k,1}}$. In other words, $\theta_0 \in \Sigma [\mathbf{x}]_1, \tilde{\theta_k} \in \Sigma [\mathbf{x}_{\Gamma_{k,1}}]_2 \subseteq \tilde{\Sigma} [\mathbf{x}]_2^2, \tilde{\sigma_k} \in \Sigma [\mathbf{x}_{\Gamma_{k,1}}]_1 \subseteq \tilde{\Sigma} [\mathbf{x}]_1^2$.
		
		Similarly, define $\Gamma_{1,1,1} = \{1,2\} \subseteq I_1$ and $\Gamma_{2,2,1} = \{5,6\} \subseteq I_2$, then the second-order sparse $(\mathbf{l}, \mathbf{q})$-sublevel relaxation reads
		$$\sup_{t \in \mathbb{R}} \{t: f(\mathbf{x}) - t = \big(\theta_{0,1} (\mathbf{x}_{I_1}) + \tilde{\theta_1} (\mathbf{x}_{\Gamma_{1,1,1}}) + \tilde{\sigma_1} (\mathbf{x}_{\Gamma_{1,1,1}}) g_1 (\mathbf{x})\big) + \big(\theta_{0,2} (\mathbf{x}_{I_2}) + \tilde{\theta_2} (\mathbf{x}_{\Gamma_{2,2,1}}) + \tilde{\sigma_2} (\mathbf{x}_{\Gamma_{2,2,1}}) g_2 (\mathbf{x})\big)\}$$
		where $\theta_{0,k}$ are degree-2 SOS polynomials in variable $\mathbf{x}_{I_k}$, $\tilde{\theta_{k}}$ are degree-4 SOS polynomials in variable $\mathbf{x}_{\Gamma_{k,k,1}}$, $\tilde{\sigma_{k}}$ are degree-2 SOS polynomials in variable $\mathbf{x}_{\Gamma_{k,k,1}}$. In other words, $\theta_{0,k} \in \Sigma [\mathbf{x}_{I_k}]_1, \tilde{\theta_k} \in \Sigma [\mathbf{x}_{\Gamma_{k,k,1}}]_2 \subseteq \tilde{\Sigma} [\mathbf{x}_{I_k}]_2^2, \tilde{\sigma_k} \in \Sigma [\mathbf{x}_{\Gamma_{k,k,1}}]_1 \subseteq \tilde{\Sigma} [\mathbf{x}_{I_k}]_1^2$.
	\end{eg}
	
	The standard Lasserre's hierarchy and many of its variants are contained in the framework of sublevel hierarchy:
	
	\begin{eg}
		\textbf{(Dense Lasserre's Relaxation \cite{lasserre2001global})} The dense version of the $d$-th order Lasserre's relaxation is the dense $(d-1)$-th order sublevel relaxation with $\mathbf{l} = [n, n, \ldots, n]$ and $\mathbf{q} = \mathbf{1}_p$, where $\mathbf{1}_p$ denotes the $p$-dimensional vector with all ones.
	\end{eg}
	
	\begin{eg}
		\textbf{(Sparse Lasserre's Relaxation \cite{lasserre2006convergent})} The sparse version of the $(d-1)$-th order Lasserre's relaxation is the sparse $d$-th order sublevel relaxation with $\mathbf{l} = [[\tau_{s}]_{s \in k(1)}; \ldots; [\tau_{s}]_{s \in k(n)}]$ and $\mathbf{q} = \mathbf{1}_{|k(1)|+\ldots + |k(n)|}$.
	\end{eg}
	
	\begin{eg}
		\textbf{(Multi-Order/Partial Relaxation)} The multi-order relaxation (used to solve the \emph{Optimal Power Flow} problem in \cite{josz2018lasserre}), also named as partial relaxation (used to solve the \emph{Max-Cut} problem in \cite{campos2020partial}), is a variant of the second-order sparse Lasserre's relaxation. We first preset a value $r$, then compute the maximal cliques in the chordal extension of the CSP graph of the POP.
		For those cliques of size larger than $r$, we consider the first-order moment matrices; for those of size smaller or equal than $r$, we consider the second-order moment matrices. Denote by $S$ the set of indices such that $\tau_k > r$ for $k \in S$, and $T$ the set of indices such that $\tau_k \le r$ for $k \in T$. Then the multi-order/partial relaxation is the second-order sublevel relaxation with $\mathbf{l} = [[0]_{s \in k(1) \cap S}, [\tau_{s}]_{s \in k(1) \cup T}; \ldots; [0]_{s \in k(n) \cap S}, [\tau_{s}]_{s \in k(n) \cup T}]$ and $\mathbf{q} = [[0]_{s \in k(1) \cap S}, [1]_{s \in k(1) \cup T}; \ldots; [0]_{s \in k(n) \cap S}, [1]_{s \in k(n) \cup T}]$.
	\end{eg}
	
	\begin{eg}
		\textbf{(Augmented Partial Relaxation)} This relaxation is the strengthened version of the partial relaxation used by the authors in \cite{campos2020partial} to solve Max-Cut problems. It is exactly the second-order sublevel relaxation restricted to Max-Cut problem.
	\end{eg}
	
	\begin{eg}
		\textbf{(Heuristic Relaxation)} The heuristic relaxation proposed by the authors in \cite{chen2020semialgebraic} to compute the upper bound of the Lipschitz constant of ReLU networks, is a variant of the second-order dense Lasserre's relaxation. The intuition is that some constraints in the POP are sparse, so let us denote by $S$ the set of their indices, while their corresponding cliques are large, thus one cannot solve the second-order relaxation of the standard sparse Lasserre's hierarchy. We then consider the dense first-order relaxation (Shor's relaxation), and choose subsets of moderate sizes (size 2 in \cite{chen2020semialgebraic}) that contain the variable sets of these sparse constraints. For other constraints with larger variable sets, let us denote by $T$ the set of their indices  and let us consider the first-order moment matrices. Then the heuristic relaxation is the second-order sublevel relaxation with $\mathbf{l} = [[0]_{i \in T}, [2]_{i \in S}]$ and $\mathbf{q} = [[0]_{i \in T}, [1]_{i \in S}]$.
	\end{eg}
	
	Summarizing the above discussion, we have the following proposition:
	\begin{prop}
		For the dense case, if $\mathbf{l} = [n, n, \ldots, n]$, then the $d$-th order $(\mathbf{l}, \mathbf{q})$-sublevel relaxation is exactly the dense $(d+1)$-th order Lasserre's relaxation.
		
		For the sparse case, if $\mathbf{l} = [[\tau_{s}]_{s \in k(1)}; \ldots; [\tau_{s}]_{s \in k(n)}]$, then the $d$-th order $(\mathbf{l}, \mathbf{q})$-sublevel relaxation is exactly the sparse $(d+1)$-th order Lasserre's relaxation.
	\end{prop}
	
	\subsection{Determining the subsets of cliques} \label{weights}
	There are different ways to determine the subsets $\Gamma_{i,j}$ (or $\Gamma_{i,k,j}$) of the sublevel relaxation described in Definition \ref{sublevel}. Generically, we are not aware of any algorithm that would  guarantee that the selected subsets are optimal at a given level of relaxation. In this section, we propose several heuristics to select the subsets. Suppose that $\{I_k\}_{1\leq k\leq r}$ is the sequence of maximal cliques in the chordal extension of the CSP graph of the sparse problem \eqref{sparse_opt} and that the level of relaxation is $l \le |I_k| =: \tau_k$. 
	We need to select the ``best'' candidate among the $\tau_k \choose l$ many subsets of size $l$. However, in practice, the number $\tau_k \choose l$ might be very large since ${\tau_k \choose l} \approx \tau_k ^l$ when $l$ is fixed. 
	
	In order to make this selection procedure tractable, we reduce the number of sample subsets to $\tau_k$. Precisely, suppose $I_k := \{i_1, i_2, \ldots, i_{\tau_k}\}$, define $I_{k,j} := \{i_j, i_{j+1}, \ldots, i_{j+l}\}$ for $j = 1, 2, \ldots, \tau_k$ and $1 \le l \le \tau_k$. By convention, $i_j = i_k$ if $j \equiv k \mod \tau_k$. Denote by $p$ the depth of the relaxation. Then we use the following heuristics to choose $p$ subsets among the candidates $I_{k,j}$. Without loss of generality, we assume that $l < \tau_k$ (otherwise one has $l \ge \tau_k$, then we only need to select one subset $I = I_k$).
	
	$\bullet$ \textbf{H1} \emph{(Random Heuristic)}. For each $i$ and clique $I_{k}$, we randomly select $p$ subsets $\Gamma_{i,k,j} \subseteq I_{k}$ for $j = 1, \ldots, p$, such that $|\Gamma_{i,k,j}| = l$ for all $j$.  
	
	$\bullet$ \textbf{H2} \emph{(Ordered Heuristic)}. For each $i$ and clique $I_{k}$, we select one after another $\Gamma_{i,k,j} = I_{k,j} \subseteq I_{k}$ for $j = 1, \ldots, p$. For $p = \tau_{k}$, we also call this heuristic the \emph{cyclic heuristic}. 
	
	The heurisics \textbf{H1} and \textbf{H2} do not depend on the problem, thus they might not fully explore the specific structure hidden in the POPs. 
	We can also try the heuristic that selects the subsets according to the value of the moments in the first-order moment relaxation (Shor's relaxation).
	
	$\bullet$ \textbf{H3} \emph{(Moment Heuristic)}. First of all, we solve the first-order sparse relaxation. For each $i$ and clique $I_k$, suppose $\mathbf{M}_k$ is the first-order moment matrix indexed by $1$ and the monomials in $\textbf{x}_{I_k}$. Denote by $\mathbf{M}_k (I_{k,j})$ the submatrix whose rows and columns are indexed by $1$ and  $\textbf{x}_{I_{k,j}}$ for $j = 1, 2, \ldots, \tau_k$. We reorder the subsets $I_{k,j}$ w.r.t. the infinity norm of the submatrices $\mathbf{M}_k (I_{k,j})$, i.e., 
	$$||\mathbf{M}_k (I_{k,1})||_{\infty} \ge ||\mathbf{M}_k (I_{k,2})||_{\infty} \ge \ldots \ge \||\mathbf{M}_k (I_{k,\tau_k})||_{\infty} \,.$$
	Then we pick the first $p$ subsets $\Gamma_{i,k,1} = I_{k,1}, \Gamma_{i,k,2} = I_{k,2}, \ldots, \Gamma_{i,k,p} = I_{k,p}$ after reordering. 
	
	In particular, for Max-Cut problem, the authors in \cite{campos2020partial} proposed the following heuristics that take the weights in the graph or the maximal cliques in the chordal graph into account. We briefly introduce the idea of these heuristics, readers can refer to \cite{campos2020partial} for details. For heuristic \textbf{H4} to \textbf{H4-6}, denote by $\mathbf{L}$ the \emph{Laplacian} matrix of the graph.
	
	$\bullet$ \textbf{H4} \emph{(Laplacian Heuristic)}. For each clique $I_k$, denote by $\mathbf{L} (I_{k,j})$ the submatrix of the moment matrix $\mathbf{M}_k$ whose rows and columns are indexed by $1$ and $\mathbf{x}_{I_{k,j}}$ for $j= 1, 2, \ldots, \tau_k$. 
	We reorder the subsets $(I_{k,j})$ w.r.t. the infinity norm of the submatrices $(\mathbf{L} (I_{k,j}))$, i.e.,
	$$||\mathbf{L} (I_{k,1})||_{\infty} \ge ||\mathbf{L} (I_{k,2})||_{\infty} \ge \ldots \ge \||\mathbf{L} (I_{k,\tau_k})||_{\infty}$$
	Then we pick the first $p$ subsets $\Gamma_{i,k,1} = I_{k,1}, \Gamma_{i,k,2} = I_{k,2}, \ldots, \Gamma_{i,k,p} = I_{k,p}$ after reordering. 
	
	$\bullet$ \textbf{H5} \emph{(Max-Repeated Heuristic)}. We select subsets contained in many maximal cliques.
	
	$\bullet$ \textbf{H6} \emph{(Min-Repeated Heuristic)}. We select subsets contained in few maximal cliques.
	
	$\bullet$ \textbf{H4-5.} We combine heuristic \textbf{H4} and \textbf{H5} to select the subsets that are not repeated in other maximal cliques and contain variables with large weights.
	
	In the spirit of the heuristic \textbf{H4-5}, we can also combine \textbf{H5} with the moment heuristic \textbf{H3}:
	
	
	$\bullet$ \textbf{H3-5.} We combine \textbf{H3} and \textbf{H5} to select the subsets that are not repeated in other maximal cliques and contain variables with large moments.
	

	%
	
	\begin{table}[H]
		\scriptsize
		\caption{Comparison of different heuristics for Max-Cut instances g\_20 and w01\_100.} \label{heuristics}
		\makebox[1 \textwidth][c]{       
			\begin{tabular}{ccccccccccccc}
				\toprule
				\multirow{2}{*}{Heuristics} & \multicolumn{2}{c}{lv=4, p=1} && \multicolumn{2}{c}{lv=4, p=2} && \multicolumn{2}{c}{lv=6, p=1} && \multicolumn{2}{c}{lv=6, p=2} & \multirow{2}{*}{Count} \\
				\cline{2-3}\cline{5-6}\cline{8-9}\cline{11-12}
				& g20 & w01 && g20 & w01 && g20 & w01 && g20 & w01 &\\
				\midrule
				\textbf{H1} & 548.4 & 725.6 && \textbf{539.0} & 720.0 && \textbf{526.6} & 709.5 && \textbf{522.0} & \textbf{700.4} & \textbf{4} \\
				\textbf{H2} & \textbf{546.4} & 728.0 && 539.9 & 721.1 && 526.9 & \textbf{705.7} && 523.1 & 701.7 & \textbf{2} \\
				H3 & 550.6 & 728.8 && 541.8 & 723.2 && 528.5 & 713.9 && 524.2 & 705.6 & 0 \\
				\textbf{H4} & 549.7 & \textbf{723.4} && 542.0 & \textbf{718.6} && 526.9 & 710.5 && 523.6 & 701.5  & \textbf{2} \\
				H5 & 553.5 & 731.0 && 543.1 & 725.8 && 529.3 & 715.6 && 525.2 & 708.4 & 0 \\
				H6 & 553.3 & 731.2 && 543.2 & 726.6 && 529.3 & 717.2 && 525.2 & 710.3 & 0 \\
				H3-5 & 550.5 & 729.5 && 541.8 & 726.6 && 528.5 & 713.8 && 524.2 & 704.8 & 0 \\
				H4-5 & 549.8 & 726.6 && 542.0 & 719.3 && 526.9 & 710.4 && 523.6 & \textbf{700.4} & 1 \\
				\bottomrule
			\end{tabular}
		}
	\end{table}
	
	There is no general guarantee that one of the heuristics always performs better than the others. In Table \ref{heuristics}, we show the upper bounds obtained by the above heuristics for two Max-Cut instances g\_20 and w01\_100 (the detail of the numerical settings and the results is referred to Section \ref{num}), for level 4, 6, and depth 1, 2, respectively. For each heuristic, we count the number of times that the heuristic performs the best. We see that, suprisingly, the random heuristic \textbf{H1} performs the best among other heuristics. The ordered heuristic \textbf{H2} and Laplacian heuristic \textbf{H4} also performs well. For the sake of simplicity, we will only consider the ordered heuristic \textbf{H2} and its variants for the forthcoming examples.

	\section{Applications of sublevel hierarchy} \label{app}
	In this section, we explicitly build different sublevel relaxations for different classes of polynomial optimization problems: Maximum Cut (Max-Cut), Maximum Clique (Max-Cliq), Mixed Integer Quadratically Constrained Programming (MIQCP) and Quadratically Constrained Quadratic Problem (QCQP). We also consider two classes of problems arising from deep learning: robustness certification and Lipschitz constant estimation of neural networks. For many deep learning applications, the targeted optimization problems are often dense or nearly-dense, due to the composition of affine maps and non-linear activation functions such as $\text{ReLU}(\mathbf{A} x) = \max \{\mathbf{A} x,0\}$. In this case, the sublevel hierarchy is indeed helpful. 
	A simple application for Lipschitz constant estimation was previously considered by the authors in \cite{chen2020semialgebraic}. 
	
	For simplicity, unless stated explicitly, we always assume that all the levels $(l_i)$ (resp. $(l_{i,k})$) and depths $(q_i)$ (resp. $(q_{i,k})$) are identical, i.e., $l_i = l, q_i = q$ for all $i$ (resp. $l_{i,k} = l, q_{i,k} = q$ for all $i,k$). 
	We say that this simplified sublevel relaxation is of \emph{level} $l$ and \emph{depth} $q$. Note that the sublevel relaxation of level 0 and depth 0 is equivalent to Shor's relaxation. 
	By convention, if $l_{i,k} \ge \tau_k$, then this sublevel $l_{i,k}$ should automatically be $\tau_k$ and the depth $p$ should be 1. For all the examples, we consider the ordered heuristic $\textbf{H2}$ or its variants to select the subsets in the sublevel relaxation.
	
	\subsection{Examples from optimization}
	The examples listed in this section are typical in optimization.
	\subsubsection*{Maximum cut (Max-Cut) problem}
	Given an undirected graph $G(V, E)$ where $V$ is a set of vertices and $E$ is a set of edges, a \emph{cut} is a partition of the vertices into two disjoint subsets. The \emph{Max-Cut} problem consists of finding a cut in a graph such that the number of edges between the two subsets is as large as possible. It can be formulated as  follows:
	\begin{align} \label{Max-Cut}
		& \max_{\mathbf{x}} \{\mathbf{x}^T \mathbf{L} \mathbf{x}: \mathbf{x} \in \{-1,1\}^n\}\,, \tag{Max-Cut}
	\end{align}
	where $\mathbf{L}$ is the Laplace matrix of the given graph of $n$ vertices, i.e., $\mathbf{L} := \diag (\mathbf{W} \mathbf{1}_n) - \mathbf{W}$ where $\mathbf{W}$ is the weight matrix of the graph. The constraints $\mathbf{x} \in \{-1,1\}^n$ are equivalent to $(x_i)^2 = 1$ for all $i$. Suppose that $(I_k)$ are the maximum cliques in the chordal extension of the given graph.
	For $i = 1, 2, \ldots, n$, denote by $k(i)$ the set of indices $s$ such that $i \in I_{s}$. For $s \in k(i)$, suppose that $I_s = \{i_1, \ldots, i_{\tau_s}\}$ so that $i_{j(i)} = i$ for $1 \le j(i) \le \tau_s$. 
	Then we select the $q$ subsets of size $l$ by order as:
	$I_{s,t} = \{i_{j(i)}, i_{j(i)+t}, \ldots, i_{j(i)+t+l-2}\}$ for $t = 1, 2, \ldots, q$. If we consider the dense sublevel hierarchy, then we directly select the subsets by order as $I_t = \{i, i+t, \ldots, i+t+l-2\}$ for $t = 1, 2, \ldots, q$.
	
	\subsubsection*{Maximum clique (Max-Cliq) problem}
	Given an undirected graph $G(V, E)$ where $V$ is a set of vertices and $E$ is a set of edges, a \emph{clique} is defined to be a set of vertices that is completely interconnected. The \emph{Max-Cliq} problem consists of determining a clique of maximum cardinality. It can be stated as a nonconvex quadratic programming problem over the unit simplex \cite{pardalos1990global} and its general formulation is:
	\begin{align} \label{Max-Cut_bis}
		& \max_{\mathbf{x}} \{\mathbf{x}^T \mathbf{A} \mathbf{x}: \sum_{i = 1}^n x_i = 1, \mathbf{x} \in [0,1]^n\}\,, \tag{Max-Cliq}
	\end{align}
	where $\mathbf{A}$ is the adjacency matrix of the given graph of $n$ vertices. The constraints $\mathbf{x} \in [0,1]^n$ are equivalent to $x^i (x^i - 1) \le 0$ for $i = 1, 2, \ldots, n$. The Max-Cliq problem is dense since we have a constraint $\sum_{i = 1}^n x_i = 1$ involving all the variables. 
	Therefore, we apply the dense sublevel hierarchy. 
	To handle the constraint $\sum_{i = 1}^n x_i = 1$, we select the $q$ subsets of size $l$ by order as $I_t = \{t, t+1, \ldots, t+l-1\}$ for $t = 1, 2, \ldots, q$. For the constraints $x_i (x_i - 1) \le 0$, we select the subsets by order as $I_t = \{i, i+t, \ldots, i+t+l-2\}$ for $t = 1, 2, \ldots, q$.
	
	\subsubsection*{Mixed integer quadratically constrained programming (MIQCP)}
	The MIQCP problem is of the following form:
	\begin{align} \label{mip}
		& \min_{\mathbf{x}} \{\mathbf{x}^T \mathbf{Q}_0 \mathbf{x} + \mathbf{b}_0^T \mathbf{x}: \mathbf{x}^T \mathbf{Q}_i \mathbf{x} + \mathbf{b}_i^T \mathbf{x} \le c_i, i = 1, \ldots, p, \nonumber \\
		& \qquad\qquad\qquad\qquad\quad \mathbf{A} \mathbf{x} = \mathbf{b}, \: \mathbf{l} \le \mathbf{x} \le \mathbf{u}, \mathbf{x}_I \in \mathbb{Z}\}\,, \tag{MIQCP}
	\end{align}
	where each $\mathbf{Q}_i$ is a symmetric matrix of size $n \times n$, $\mathbf{A}$ is a matrix of size $n\times n$, $\mathbf{b}, \mathbf{b}_i, \mathbf{l}, \mathbf{u}$ are $n$-dimensional vectors, and each $c_i$ is a real number. 
	The constraints $\mathbf{x}^T \mathbf{Q}_i \mathbf{x} + \mathbf{b}_i^T \mathbf{x} \le c_i$ are called \emph{quadratic} constraints, the constraints $\mathbf{A} \mathbf{x} = \mathbf{b}$ are called \emph{linear} constraints. The constraints $\mathbf{l} \le \mathbf{x} \le \mathbf{u}$ and $x_I \in \mathbb{Z}$ bound the variables and restrict some of them to be integers. In our benchmarks, we only consider the  case where $\mathbf{x} \in \{0, 1\}^n$, which is also equivalent to $x_i (x_i - 1) = 0$ for $i = 1, 2, \ldots, n$. 
	If we only have bound constraints, then we use the same ordered heuristic as for the Max-Cut problem to select the subsets. 
	If in addition we also have quadratic constraints or linear constraints, then the problem is dense and therefore we consider the dense sublevel hierarchy. For quadratic constraints, we don't apply the sublevel relaxation to them, i.e., $l=q=0$. However, if $\mathbf{Q}_i$ equals the identity matrix, then we use the same heuristic as the linear constraints: we select the subsets by order as $I_t = \{t, t+1, \ldots, t+l-1\}$ for $t = 1, 2, \ldots, q$.

	\subsubsection*{Quadratically constrained quadratic problems (QCQP)}
	A QCQP can be cast as follows:
	\begin{align} \label{qcqp}
		& \min_{\mathbf{x}} \{\mathbf{x}^T \mathbf{Q}_0 \mathbf{x} + \mathbf{b}_0^T \mathbf{x}: \mathbf{x}^T \mathbf{Q}_i \mathbf{x} + \mathbf{b}_i^T \mathbf{x} \le c_i, i = 1, \ldots, p, \nonumber \\
		& \qquad\qquad\qquad\qquad\quad \mathbf{A} \mathbf{x} = \mathbf{b}, \: \mathbf{l} \le \mathbf{x} \le \mathbf{u}\}\,, \tag{QCQP}
	\end{align}
	where each $\mathbf{Q}_i$ is a symmetric matrix of size $n \times n$, $\mathbf{A}$ is a matrix of size $n\times n$, $\mathbf{b}, \mathbf{b}_i, \mathbf{l}, \mathbf{u}$ are $n$-dimensional vectors, and each $c_i$ is a real number. 
	This is very similar to the MIQCP except that we drop out the integer constraints. Therefore, we use the same strategy to select the subsets in the sublevel relaxation.
	
	\subsection{Examples from deep learning}
	The following examples are picked from the recent deep learning topics.
	\subsubsection*{Upper bounds of lipschitz constants of deep neural networks \cite{chen2020semialgebraic}}
	We only consider the 1-hidden layer neural network with ReLU activation function, the upper bound of whose Lipschitz constant results in a QCQP as follows:
	\begin{align} \label{lcep}
		& \max_{\mathbf{x},\mathbf{u},\mathbf{t}} \{\mathbf{t}^T \mathbf{A}^T \diag (\mathbf{u}) \mathbf{c}: \mathbf{u} (\mathbf{u} - 1) = 0, (\mathbf{u} - 1/2) (\mathbf{A} \mathbf{x} + \mathbf{b}) \ge 0; \nonumber\\
		& \qquad\qquad\qquad\qquad\quad \mathbf{t}^2 \le 1, (\mathbf{x}- \bar{\mathbf{x}} + \varepsilon) (\mathbf{x} - \bar{\mathbf{x}} - \varepsilon) \le 0\,.\} \tag{Lip}
	\end{align}
	where $\mathbf{A}$ is a matrix of size $p_2 \times p_1$, $\bar{\mathbf{x}}$ is a $p_1$-dimensional vector, $\mathbf{b}$, $\mathbf{c}$ are $p_2$-dimensional vectors, and $\epsilon$ is a positive real number. When $\epsilon = 10$ (resp. $\epsilon = 0.1$), we compute the upper bounds of the \emph{global} (resp. \emph{local}) Lipschitz constant of the neural network. Assume the matrix $\mathbf{A}$ is dense, then the maximal cliques in the chordal extension of \eqref{lcep} are $I = \{x_1, \ldots, x_{p_1}; u_1, \ldots, u_{p_2}\}$ and $I_k = \{u_1, \ldots, u_{p_2}, t_k\}$ for $k = 1, \ldots, p_1$.
	Therefore, we consider the sparse sublevel relaxation. For the constraints $t_k^2 \le 1$, we choose the subsets by order as $I_{k,i} = \{u_{i}, \ldots, u_{i+l-2}; t_k\}$ for $i = 1, \ldots, q$. 
	For the constraints $(x_i- \bar{x}_k + \varepsilon) (x - \bar{x}_k - \varepsilon) \le 0$, we choose the subsets by order as $I_i = \{x_k, x_{k+i}, \ldots, x_{k+i+l/2-2}; u_i, \ldots, u_{i+l/2-1}\}$ for $i = 1, \ldots, q$. 
	For the constraints $u_j(u_j -1) = 0$ and $(u_j -1/2) (\mathbf{A}^{j,:} \mathbf{x} + b_j) \ge 0$, we choose the subsets by order as $I_i = \{x_i, \ldots, x_{i+l/2-1}; u_j, u_{j+i}, \ldots, u_{j+i+l/2-2}\}$ for $i = 1, \ldots, q$.
	
	\subsection*{Robustness certification of deep neural networks \cite{Raghuathan18}}
	We also consider the 1-hidden layer neural network with ReLU activation function. 
	Then the robustness certification problem can be formulated as a QCQP as follows:
	\begin{align} \label{cert}
		& \max_{\mathbf{x},\mathbf{u}} \{\mathbf{c}^T \mathbf{u}: \mathbf{u} (\mathbf{u} - \mathbf{A} \mathbf{x} - \mathbf{b}) = 0, \: \mathbf{u} \ge \mathbf{A} \mathbf{x} + \mathbf{b}, \: \mathbf{u} \ge 0, \nonumber \\
		& \qquad\qquad\quad (\mathbf{x}- \bar{\mathbf{x}} + \varepsilon) (\mathbf{x} - \bar{\mathbf{x}} - \varepsilon) \le 0\,.\} \tag{Cert}
	\end{align}
	where $\mathbf{A}$ is a matrix of size $p_2 \times p_1$, $\bar{\mathbf{x}}$ is a $p_1$-dimensional vector, $\mathbf{b}$, $\mathbf{c}$ are $p_2$-dimensional vectors, and $\epsilon$ is a positive real number. Assume the matrix $\mathbf{A}$ is dense, then the maximal cliques in the chordal extension of \eqref{cert} are $I_k = \{x_1, \ldots, x_{p_1}; u_k\}$ for $k = 1, \ldots, p_2$.
	Similarly to the Lipschitz problem \eqref{lcep}, we consider the sparse sublevel relaxation. For all the constraints, we choose the subsets by order as $I_{k,i} = \{x_i, \ldots, x_{i+l-2}; u_k\}$ for $i = 1, \ldots, q$.
	
	
	\section{Numerical results} \label{num}
	In this section, we apply the sublevel relaxation to different type of POPs both in optimization and deep learning, as discussed in the previous section. Most of the instances in optimization are taken from the Biq-Mac library \cite{rendl2007branch} and the QPLIB library \cite{furini2019qplib}, others are generated randomly. 
	We calculate the \emph{ratio of improvements (RI)} of each sublevel relaxation, compared with Shor's relaxation, namely $\text{RI} = \frac{\text{Shor} - \text{sublevel}}{\text{Shor} - \text{solution}} \times 100\%$. We also compute the \emph{relative gap (RG)} between the sublevel relaxation and the optimal solution, given by $\text{RG} = \frac{\text{sublevel} - \text{solution}}{|\text{solution}|} \times 100\%$. For each instance, we only show the ratio of improvements and relative gap corresponding to the results of the last sublevel relaxation. The larger the ratio of improvements or the smaller the relative gap, the better the bounds. If the optimal solution is not known so far, it is replaced by the (best-known) valid upper bounds (UB) or lower bounds (LB). We implement all the programs on Julia, and use Mosek as back-end to solve SDP relaxations. The running time (with second as unit) displayed in all tables refers to the time spent by Mosek to solve the SDP relaxation. 
	All experiments are performed with an Intel 8-Core i7-8665U CPU @ 1.90GHz Ubuntu 18.04.5 LTS, 32GB RAM.
	
	%

	\subsection{Examples from optimization}
	\subsubsection*{Max-Cut instances}\label{max}
	The following classes of problems and their solutions are from the Biq-Mac library. 
	For each class of problem, we choose the first instance, i.e., $i = 0$, and drop the suffix ``.i'' in Table \ref{maxcut}:
	
	$\bullet$ g05\_$n.i$, unweighted graphs with edge probability 0.5, $n = 60, 80, 100$.
	
	$\bullet$ pm1s\_$n.i$, pm1d\_$n.i$, weighted graph with edge weights chosen uniformly from $\{-1,0,1\}$ and density 10\% and 99\% respectively, $n= 80,100$.
	
	$\bullet$ w$d\_n.i$, pw$d\_n.i$, graph with integer edge weights chosen from $[-10,10]$ and $[0,10]$ respectively, density $d = 0.1, 0.5, 0.9$, $n = 100$.
	
	The instances named g\_$n$ and the corresponding upper bounds are from the CS-TSSOS paper \cite{wang2020cs}.
	
	The instances named G$n$ are from the G-set library by Y.Y. Ye \footnote{\url{http://web.stanford.edu/~yyye/yyye/Gset/}}, and their best-known solutions are taken from \cite{kochenberger2013solving}.
	
	In Table \ref{maxcut-summary}, we give a summary of basic information and the graph structure of each instance: \emph{nVar} denotes the number of variables, \emph{Density} denotes the percentage of non-zero elements in the adjacency matrix, \emph{nCliques} denotes the number of cliques in the chordal extension, \emph{MaxClique} denotes the maximum size of the cliques, \emph{MinClique} denotes the minimum size of the cliques.

	\begin{table}[H]
		\scriptsize
		\caption{Summary of the basic information and graph structure of the Max-Cut instances.} \label{maxcut-summary}
		\makebox[1 \textwidth][c]{       
			\begin{tabular}{|c|c|c|c|c|c|}
				\hline
				& nVar & Density & nCliques & MaxClique & MinClique \\
				\hline
				g05\_60 & 60 & 50\% & 11 & 50 & 19 \\
				\hline
				g05\_80 & 80 & 50\% & 12 & 69 & 28 \\
				\hline
				g05\_100 & 100 & 50\% & 13 & 88 & 37 \\
				\hline 
				pm1d\_80 & 80 & 99\% & 2 & 79 & 76 \\
				\hline
				pm1d\_100 & 100 & 99\% & 2 & 99 & 95 \\
				\hline
				pm1s\_80 & 80 & 10\% & 44 & 37 & 4 \\
				\hline
				pm1s\_100 & 100 & 10\% & 47 & 54 & 4 \\
				\hline
				pw01\_100 & 100 & 10\% & 47 & 54 & 4 \\
				\hline
				pw05\_100 & 100 & 50\% & 12 & 89 & 40 \\
				\hline
				pw09\_100 & 100 & 90\% & 4 & 97 & 83 \\
				\hline
				w01\_100 & 100 & 10\% & 47 & 54 & 4\\
				\hline
				w05\_100 & 100 & 50\% & 12 & 89 & 40 \\
				\hline
				w09\_100 & 100 & 90\% & 4 & 97 & 83 \\
				\hline
				g\_20 & 505 & 1.6\% & 369 & 15 & 1\\
				\hline
				g\_40 & 1005 & 0.68\% & 756 & 15 & 1\\
				\hline
				g\_60 & 1505 & 0.43\% & 756 & 15 & 1\\
				\hline
				g\_80 & 2005 & 0.30\% & 1556 & 15 & 1\\
				\hline
				g\_100 & 2505 & 0.23\% & 1930 & 16 & 1\\
				\hline
				g\_120 & 3005 & 0.19\% & 2383 & 15 & 1\\
				\hline
				g\_140 & 3505 & 0.16\% & 2762 & 15 & 1\\
				\hline
				g\_160 & 4005 & 0.13\% & 3131 & 15 & 1\\
				\hline
				g\_180 & 4505 & 0.12\% & 3429 & 15 & 1\\
				\hline
				g\_200 & 5005 & 0.11\% & 3886 & 15 & 1\\
				\hline
				G11 & 800 & 0.25\% & 598 & 24 & 5 \\
				\hline
				G12 & 800 & 0.25\% & 598 & 48 & 5 \\
				\hline
				G13 & 800 & 0.25\% & 598 & 90 & 5 \\
				\hline
				G32 & 2000 & 0.1\% & 1498 & 76 & 5 \\
				\hline
				G33 & 2000 & 0.1\% & 1498 & 99 & 5 \\
				\hline
				G34 & 2000 & 0.1\% & 1498 & 141 & 5 \\
				\hline
			\end{tabular}
		}
	\end{table}
	
	In Table \ref{maxcut}, we display the upper bounds and running times corresponding to the sublevel relaxations of depth 1, and level 0, 4, 6, 8, respectively. Notice that the authors in \cite{campos2020partial} use the partial relaxation to compute upper bounds for instances g\_20 to g\_200. 
	The sublevel relaxation we consider here is actually what they call the augmented partial relaxation, which is a strengthened relaxation based on partial relaxation. 
	From the ratio of improvement, we see that the more sparse structure the graph has, the better the sublevel relaxation performs. 
	Notice that if we obtain better upper/lower bounds than the current best-known bounds, the ratio of improvements will be larger than 100\% and the relative gap will become nagative. Particularly, our method provides better bounds for all the instances g\_$n$ in the CS-TSSOS paper \cite{wang2020cs}, and computes upper bounds very close to the best-known solution for the instances G$n$ in G-set.

	Moreover, if the number of variables is of moderate size, the dense sublevel relaxation might performs faster than the sparse one. For example, the instance g05\_100 has 13 maximal cliques with maximum size 88 and minimum size 37. The sparse sublevel relaxation consists of 13 first-order moment matrices of size from 37 to 88. However, the dense version only consists of 1 first-order moment matrix of size 100.
	In fact, the dense sublevel relaxation gives an upper bound of 1463.5 at level 0 in 10 seconds, yielding the same bound as the sparse case at level 0 but with much less computing time, and 1458.1 at level 8 in 178.1 seconds, providing better bounds than the sparse case at level 6, with less computing time.

	\begin{table}[H]
		\scriptsize
		\caption{Results obtained with sublevel relaxations of Max-Cut problems.} \label{maxcut}
		\makebox[1 \textwidth][c]{       
			\begin{tabular}{|c|c|c|c|ccccc|cccc|}
				\hline
				& \multirow{2}{*}{Sol./UB} & \multirow{2}{*}{nVar} & \multirow{2}{*}{Density} & \multicolumn{9}{c|}{Sublevel relaxation, $l$ = 0/4/6/8, $q$ = 1 (level 0 = Shor)} \\
				\cline{5-13}
				& &&& \multicolumn{5}{c|}{upper bounds (RI, RG)} & \multicolumn{4}{c|}{solving time (s)} \\
				\hline
				g05\_60 & 536 & 60 & 50\% & 550.1 & 548.1 & 546.0 & 544.6 & (39.0\%, 1.6\%) & 4.5 & 10.6 & 17.6 & 65.7 \\
				\hline
				g05\_80 & 929 & 80 & 50\% & 950.9 & 949.0 & 946.6 & 944.6 & (28.8\%, 1.7\%) & 33.8 & 56.2 & 61.8 & 137.4 \\
				\hline
				g05\_100 & 1430 & 100 & 50\% & 1463.5 & 1462.0 & 1459.2 & 1456.8 & (20.0\%, 1.9\%) & 138.7 & 303.7 & 328.7 & 460.3  \\
				\hline
				pm1d\_80 & 227 & 80 & 99\% & 270.0 & 265.9 & 262.0 & 258.8 & (26.0\%, 14.0\%) & 15.0 & 29.4 & 39.2 & 128.1  \\
				\hline
				pm1d\_100 & 340 & 100 & 99\% & 405.4 & 402.2 & 397.9 & 393.7 & (19.0\%, 15.8\%) & 47.6 & 69.4 & 110.2 & 225.1   \\
				\hline
				pm1s\_80 & 79 & 80 & 10\% & 90.3 & 86.7 & 83.6 & \textbf{82.8} & (\textbf{66.4\%}, 4.8\%) & 1.4 & 4.9 & 13.4 & 37.7   \\
				\hline
				pm1s\_100 & 127 & 100 & 10\% & 143.2 & 141.4 & 137.6 & \textbf{135.3} & (\textbf{48.8\%}, 6.5\%)& 11.1 & 24.3 & 28.6 & 180.3  \\
				\hline
				pw01\_100 & 2019 & 100 & 10\% & 2125.4 & 2107.8 & 2088.1 & \textbf{2075.0} & (\textbf{47.4\%}, 2.8\%) & 13.0 & 20.5 & 29.7 & 285.8  \\
				\hline
				pw05\_100 & 8190 & 100 & 50\% & 8427.7 & 8416.6 & 8403.6 & 8388.1 & (16.7\%, 2.4\%) & 136.8 & 223.0 & 272.9 & 400.3  \\
				\hline
				pw09\_100 & 13585 & 100 & 90\% & 13806.0 & 13797.1 & 13781.1 & 13766.5 & (17.9\%, 1.3\%) & 141.6 & 218.4 & 268.7 & 442.4 \\
				\hline
				w01\_100 & 651 & 100 & 10\% & 740.9 & 728.3 & 710.3 & \textbf{696.2} & (\textbf{49.7\%}, 6.9\%) & 10.5 & 22.4 & 35.0 & 224.7  \\
				\hline
				w05\_100 & 1646 & 100 & 50\% & 1918.0 & 1902.6 & 1885.5 & 1869.7 & (17.8\%, 13.6\%) & 138.1 & 265.8 & 272.2 & 403.2  \\
				\hline
				w09\_100 & 2121 & 100 & 90\% & 2500.3 & 2478.2 & 2447.3 & 2422.8 & (20.4\%, 14.2\%) & 124.3 & 255.0 & 280.8 & 451.7  \\
				\hline
				g\_20 & 537.4 & 505 & 1.6\% & 570.8 & 547.1 & 526.7 & \textbf{513.4} & (\textbf{171.9\%}, -4.5\%) & 0.7 & 15.1 & 46.1 & 102.2   \\
				\hline
				g\_40 & 992.2 & 1005 & 0.68\% & 1032.6 & 982.4 & 950.8 & \textbf{927.6} & (\textbf{260.0\%}, -6.5\%) & 1.2 & 18.6 & 47.9 & 102.5   \\
				\hline
				g\_60 & 1387.2 & 1505 & 0.43\% & 1439.9 & 1368.4 & 1317.8 & \textbf{1281.9} & (\textbf{300.4\%}, -7.6\%) & 2.8 & 26.0 & 74.7 & 431.1   \\
				\hline
				g\_80 & 1838.1 & 2005 & 0.3\% & 1899.2 & 1803.8 & 1744.9 & \textbf{1698.8} & (\textbf{328.0\%}, -7.6\%) & 6.0 & 23.8 & 76.0 & 290.7   \\
				\hline
				g\_100 & 2328.3 & 2505 & 0.23\% & 2398.7 & 2282.9 & 2205.1 & \textbf{2149.3} & (\textbf{354.3\%}, -7.7\%) & 3.4 & 30.1 & 117.4 & 428.6   \\
				\hline
				g\_120 & 2655.4 & 3005 & 0.19\% & 2731.7 & 2588.5 & 2507.3 & \textbf{2439.8} & (\textbf{382.6\%}, -8.1\%) & 3.8 & 33.3 & 113.2 & 434.5   \\
				\hline
				g\_140 & 3027.2 & 3505 & 0.16\% & 3115.8 & 2947.9 & 2856.5 & \textbf{2782.6} & (\textbf{376.1\%}, -8.1\%) & 3.8 & 46.3 & 138.4 & 522.1   \\
				\hline
				g\_160 & 3589.0 & 4005 & 0.13\% & 3670.7 & 3487.1 & 3380.7 & \textbf{3310.9} & (\textbf{440.4\%}, -7.7\%) & 8.2 & 56.5 & 198.2 & 506.6   \\
				\hline
				g\_180 & 3953.1 & 4505 & 0.12\% & 4054.7 & 3855.9 & 3736.9 & \textbf{3653.5} & (\textbf{394.9\%}, -7.6\%) & 8.8 & 51.5 & 277.0 & 693.4   \\
				\hline
				g\_200 & 4472.3 & 5005 & 0.11\% & 4584.6 & 4353.3 & 4228.1 & \textbf{4132.2} & (\textbf{402.8\%}, -7.6\%) & 5.4 & 52.7 & 203.2 & 839.2   \\
				\hline
				G11 & 564 & 800 & 0.25\% & 629.2 & 581.3 & 564.6 & \textbf{564.6} & (\textbf{99.1\%}, 0.1\%) & 4.0 & 15.8 & 32.6 & 36.5  \\
				\hline
				G12 & 556 & 800 & 0.25\% & 623.9 & 572.5 & 559.6 & \textbf{559.6} & (\textbf{94.7\%}, 0.6\%) & 17.8 & 57.8 & 54.3 & 51.9  \\
				\hline
				G13 & 580 & 800 & 0.25\% & 647.1 & 594.2 & 585.1 & \textbf{584.1} & (\textbf{93.9\%}, 0.7\%) & 159.2 & 241.7 & 340.2 & 321.6  \\
				\hline
				G32 & 1398 & 2000 & 0.1\% & 1567.6 & 1433.4 & 1415.9 & \textbf{1415.9} & (\textbf{89.4\%}, 1.3\%) & 622.0 & 736.3 & 630.8 & 628.0  \\
				\hline
				G33 & 1376 & 2000 & 0.1\% & 1544.3 & 1415.3 & 1392.7 & \textbf{1387.4} & (\textbf{93.2\%}, 0.8\%) & 1956.6 & 2115.8 & 1221.5 & 1486.8  \\
				\hline
				G34 & 1372 & 2000 & 0.1\% & 1546.7 & 1407.9 & 1388.2 & \textbf{1388.2} & (\textbf{90.7\%}, 1.2\%) & 3613.5 & 6580.9 & 6327.9 & 6147.4  \\
				\hline
			\end{tabular}
		}
	\end{table}

	\subsubsection*{MIQCP instances}
	The following classes of problems and their solutions are from the Biq-Mac library, where there are neither quadratic constraints $\mathbf{x}^T \mathbf{Q}_i \mathbf{x} + \mathbf{b}_i^T \mathbf{x} \le c_i$ nor linear constraints $\mathbf{A} \mathbf{x} = \mathbf{b}$. We only have integer bound constraints $\mathbf{x} \in \{0, 1\}^n$.
	
	$\bullet$ bqp$n$-$i$, with 10\% density. All the coefficients have uniformly chosen integer values in $[-100,100]$, $n = 50, 100, 250, 500$.
	
	$\bullet$ gka$i$a, with dimensions in $[30,100]$ and densities in $[0.0625, 0.5]$. The diagonal coefficients lie in $[-100,100]$ and the off-diagonal coefficients belong to $[-100,100]$.
	
	$\bullet$ gka$i$b, with dimensions in $[20, 125]$ and density 1. The diagonal coefficients lie in $[-63,0]$ and the off-diagonal coefficients belong to $[0,100]$.
	
	$\bullet$ gka$i$c, dimensions in $[40, 100]$ and densities in $[0.1, 0.8]$. Diagonal coefficients in $[-100,100]$, off-diagonal coefficients in $[-50,50]$.
	
	$\bullet$ gka$i$d, with dimension 100 and densities in $[0.1, 1]$. The diagonal coefficients lie in $[-75,75]$ and the off-diagonal coefficients belong to $[-50,50]$.
	
	We also select some instances and their solutions from the QPLIB library with ID 0032, 0067, 0633, 2512, 3762, 5935 and 5944, in which we have additional linear constraints $\mathbf{A} \mathbf{x} = \mathbf{b}$. For the instance 0032, there are 50 continuous variables and 50 integer variables. For the two instances 5935 and 5944, we maximize the objective, the others are minimization problems.

	Similarly to the Max-Cut instances, Table \ref{miqcp_summary} summarizes the basic information and cliques structure of each instance. Table \ref{qp_summary} is a summary of basic information and the number of quadratic, linear, bound constraints of the instances from the QPLIB library. 
	
	\begin{table}[H]
		\scriptsize
		\caption{Summary of the basic information and sparse structure of the MIQCP instances.} \label{miqcp_summary}
		\makebox[1 \textwidth][c]{       
			\begin{tabular}{|c|c|c|c|c|c|c|c|c|}
				\hline
				& nVar & Density & nCliques & MaxClique & MinClique & nQuad & nLin & nBound \\
				\hline
				bqp50-1 & 50 & 10\% & 36 & 15 & 3 & 0 & 0 & 50 \\
				\hline
				bqp100-1 & 100 & 10\% & 52 & 49 & 4 & 0 & 0 & 100 \\
				\hline
				gka1a & 50 & 10\% & 36 & 15 & 1 & 0 & 0 & 50 \\
				\hline
				gka2a & 60 & 10\% & 41 & 20 & 3 & 0 & 0 & 60 \\
				\hline
				gka3a & 70 & 10\% & 44 & 27 & 3 & 0 & 0 & 70 \\
				\hline
				gka4a & 80 & 10\% & 48 & 33 & 4 & 0 & 0 & 80 \\
				\hline
				gka5a & 50 & 20\% & 25 & 26 & 4 & 0 & 0 & 50 \\
				\hline
				gka6a & 30 & 40\% & 11 & 20 & 7 & 0 & 0 & 30 \\
				\hline
				gka7a & 30 & 50\% & 10 & 21 & 10 & 0 & 0 & 30 \\
				\hline
				gka8a & 100 & 62.5\% & 64 & 37 & 2 & 0 & 0 & 100 \\
				\hline
				gka1b & 20 & 100\% & 2 & 19 & 19 & 0 & 0 & 20 \\
				\hline
				gka2b & 30 & 100\% & 2 & 29 & 29 & 0 & 0 & 30 \\
				\hline
				gka3b & 40 & 100\% & 2 & 39 & 38 & 0 & 0 & 40 \\
				\hline
				gka4b & 50 & 100\% & 2 & 49 & 47 & 0 & 0 & 50 \\
				\hline
				gka5b & 60 & 100\% & 2 & 59 & 56 & 0 & 0 & 60 \\
				\hline
				gka6b & 70 & 100\% & 2 & 69 & 67 & 0 & 0 & 70 \\
				\hline
				gka7b & 80 & 100\% & 2 & 79 & 77 & 0 & 0 & 80 \\
				\hline
				gka8b & 90 & 100\% & 2 & 89 & 87 & 0 & 0 & 90 \\
				\hline
				gka9b & 100 & 100\% & 2 & 99 & 97 & 0 & 0 & 100 \\
				\hline
				gka10b & 125 & 100\% & 2 & 124 & 124 & 0 & 0 & 125\\
				\hline
				gka1c & 40 & 80\% & 4 & 37 & 25 & 0 & 0 & 40 \\
				\hline
				gka2c & 50 & 60\% & 6 & 45 & 26 & 0 & 0 & 50 \\
				\hline
				gka3c & 60 & 40\% & 14 & 47 & 17 & 0 & 0 & 60 \\
				\hline
				gka4c & 70 & 30\% & 22 & 49 & 12 & 0 & 0 & 70 \\
				\hline
				gka5c & 80 & 20\% & 27 & 54 & 11 & 0 & 0 & 80 \\
				\hline
				gka6c & 90 & 10\% & 46 & 45 & 4 & 0 & 0 & 90 \\
				\hline
				gka7c & 100 & 10\% & 51 & 50 & 3 & 0 & 0 & 100 \\
				\hline
				gka1d & 100 & 10\% & 50 & 51 & 4 & 0 & 0 & 100 \\
				\hline
				gka2d & 100 & 20\% & 30 & 71 & 11 & 0 & 0 & 100 \\
				\hline
				gka3d & 100 & 30\% & 23 & 78 & 18 & 0 & 0 & 100 \\
				\hline
				gka4d & 100 & 40\% & 15 & 86 & 31 & 0 & 0 & 100 \\
				\hline
				gka5d & 100 & 50\% & 13 & 88 & 36 & 0 & 0 & 100 \\
				\hline
				gka6d & 100 & 60\% & 10 & 91 & 47 & 0 & 0 & 100 \\
				\hline
				gka7d & 100 & 70\% & 7 & 94 & 57 & 0 & 0 & 100 \\
				\hline
				gka8d & 100 & 80\% & 6 & 95 & 68 & 0 & 0 & 100 \\
				\hline
				gka9d & 100 & 90\% & 5 & 96 & 79 & 0 & 0 & 100 \\
				\hline
				gka10d & 100 & 100\% & 2 & 99 & 95 & 0 & 0 & 100 \\
				\hline 
			\end{tabular}
		}
	\end{table}
	
	\begin{table}[H]
		\scriptsize
		\caption{Summary of the basic information and constraint structure of the MIQCP instances from QPLIB library.} \label{qp_summary}
		\makebox[1 \textwidth][c]{       
			\begin{tabular}{|c|c|c|c|c|c|c|}
				\hline
				& nVar & Density & nQuad & nLin & nBound \\
				\hline
				qplib0032 & 100 & 89\% & 0 & 52 & 100 \\
				\hline
				qplib0067 & 80 & 89\% & 0 & 1 & 80 \\
				\hline
				qplib0633 & 75 & 99\% & 0 & 1 & 75 \\
				\hline
				qplib2512 & 100 & 28\% & 0 & 20 & 100 \\
				\hline
				qplib3762 & 90 & 28\% & 0 & 480 & 90 \\
				\hline
				qplib5935 & 100 & 28\% & 0 & 1237 & 100 \\
				\hline
				qplib5944 & 100 & 28\% & 0 & 2475 & 100 \\
				\hline
			\end{tabular}
		}
	\end{table}

	In Table \ref{miqcp}, we show the lower bounds and running time obtained by solving the sublevel relaxations with depth 1 and level 0, 4, 6, 8, respectively. 
	We see that when the problem has a good sparsity structure or is of low dimension, the sublevel relaxation performs very well and provides the exact solution, in particular for the two instances gka2a and gka7a. For dense problems, we are not able to find the exact solution, but still have improvements between 20\% and 40\%  compared to Shor's relaxation. Notice that for the instances gka1b to gka10b, even though we have an improvement ratio ranging from 24.0\% to 77.9\%, the relative gap is very high, varying from 38.2\% to 947.2\%. This means that these problems themselves are very hard to solve, so that the gap between the results of Shor's relaxation and the exact optimal solution is very large. Even though  the sublevel relaxation yields substantial improvement compared to Shor's relaxation, it's still far away from the true optimum.

	\begin{table}[H]
		\scriptsize
		\caption{Results obtained with sublevel relaxations of MIQCP problems.} \label{miqcp}
		\makebox[1 \textwidth][c]{       
			\begin{tabular}{|c|c|c|c|ccccc|cccc|}
				\hline
				& \multirow{2}{*}{Sol.} & \multirow{2}{*}{nVar} & \multirow{2}{*}{Density} & \multicolumn{9}{c|}{Sublevel relaxation, $l$ = 0/4/6/8, $q$ = 1 (level 0 = Shor)} \\
				\cline{5-13}
				& &&& \multicolumn{5}{c|}{lower/upper bounds (RI, RG)} & \multicolumn{4}{c|}{solving time (s)} \\
				\hline
				bqp50-1 & -2098 & 50 & 10\% & -2345.5 & -2136.3 & -2116.3 & \textbf{-2105.4} & (\textbf{97.0\%}, 0.4\%) & 0.1 & 0.5 & 1.3 & 3.4   \\
				\hline
				bqp100-1 & -7970 & 100 & 10\% & -8721.1 & -8358.2 & -8215.1 & \textbf{-8101.8} & (\textbf{82.5\%}, 1.7\%) & 8.8 & 16.7 & 21.8 & 87.9   \\
				\hline
				gka1a & -3414 & 50 & 10\% & -3623.3 & -3453.2 & -3432.6 & \textbf{-3428.5} & (\textbf{93.1\%}, 0.4\%) & 0.1 & 0.6 & 0.9 & 1.8   \\
				\hline
				gka2a & -6063 & 60 & 10\% & -6204.3 & -6076.3 & -6063.0 & \textbf{-6063.0} & (\textbf{100\%}, 0\%) & 0.3 & 1.0 & 4.6 & 8.5   \\
				\hline
				gka3a & -6037 & 70 & 10\% & -6546.2 & -6291.5 & -6182.6 & \textbf{-6106.3} & (\textbf{86.4\%}, 1.1\%) & 0.7 & 1.6 & 6.1 & 31.0   \\
				\hline
				gka4a & -8598 & 80 & 10\% & -8935.1 & -8767.3 & -8713.7 & \textbf{-8676.0} & (\textbf{76.9\%}, 0.9\%) & 2.1 & 3.4 & 10.1 & 30.0   \\
				\hline
				gka5a & -5737 & 50 & 20\% & -5979.9 & -5789.9 & -5760.3 & \textbf{-5750.0} & (\textbf{94.6\%}, 0.2\%) & 0.7 & 1.4 & 6.2 & 31.0   \\
				\hline
				gka6a & -3980 & 30 & 40\% & -4190.2 & -4008.9 & -3986.0 & \textbf{-3982.5} & (\textbf{98.8\%}, 0.1\%) & 0.2 & 0.6 & 3.9 & 23.6   \\
				\hline
				gka7a & -4541 & 30 & 50\% & -4696.6 & -4566.8 & -4541.1 & \textbf{-4541.1} & (\textbf{100\%}, 0\%) & 0.3 & 0.8 & 4.9 & 23.1   \\
				\hline
				gka8a & -11109 & 100 & 62.5\% & -11283.8 & -11148.0 & -11124.8 & \textbf{-11114.0} & (\textbf{97.1\%}, 0.05\%) & 2.3 & 2.7 & 7.5 & 19.5   \\
				\hline
				gka1b & -133 & 20 & 100\% & -362.9 & -295.1 & -253.6 & \textbf{-183.8} & (\textbf{77.9\%}, 38.2\%) & 0.1 & 0.5 & 2.4 & 25.0   \\
				\hline
				gka2b & -121 & 30 & 100\% & -505.7 & -425.3 & -325.4 & \textbf{-282.5} & (\textbf{58.0\%}, 133.5\%) & 0.2 & 0.7 & 4.0 & 29.9   \\
				\hline
				gka3b & -118 & 40 & 100\% & -718.0 & -535.6 & -483.4 & -437.7 & (\textbf{46.7\%}, 270.9\%) & 0.7 & 1.4 & 6.5 & 45.9   \\
				\hline
				gka4b & -129 & 50 & 100\% & -809.8 & -670.9 & -614.2 & -571.5 & (35.0\%, 343.0\%) & 1.9 & 3.3 & 14.3 & 65.2   \\
				\hline
				gka5b & -150 & 60 & 100\% & -1034.8 & -820.9 & -736.8 & -705.5 & (37.2\%, 370.3\%) & 3.2 & 8.4 & 15.5 & 76.1   \\
				\hline
				gka6b & -146 & 70 & 100\% & -1279.0 & -972.2 & -894.8 & -833.5 & (39.3\%, 470.9\%) & 9.1 & 11.6 & 26.6 & 86.5   \\
				\hline
				gka7b & -160 & 80 & 100\% & -1362.5 & -1138.1 & -1031.0 & -982.6 & (31.6\%, 514.1\%) & 26.1 & 31.2 & 50.8 & 136.1   \\
				\hline
				gka8b & -145 & 90 & 100\% & -1479.1 & -1269.8 & -1190.2 & -1120.9 & (26.8\%, 673.0\%) & 40.5 & 60.1 & 102.3 & 187.0   \\
				\hline
				gka9b & -137 & 100 & 100\% & -1663.6 & -1385.4 & -1298.9 & -1212.6 & (29.5\%, 785.1\%) & 65.9 & 92.3 & 111.2 & 256.3   \\
				\hline
				gka10b & -154 & 125 & 100\% & -2073.1 & -1782.1 & -1707.1 & -1612.7 & (24.0\%, 947.2\%) & 285.8 & 413.3 & 452.2 & 700.9   \\
				\hline
				gka1c & -5058 & 40 & 80\% & -5161.1 & -5102.9 & -5077.9 & \textbf{-5073.7} & (\textbf{84.8\%}, 0.3\%) & 0.8 & 1.6 & 5.3 & 41.9   \\
				\hline
				gka2c & -6213 & 50 & 60\% & -6392.6 & -6291.3 & -6263.1 & \textbf{-6246.2} & (\textbf{81.5\%}, 0.5\%) & 1.9 & 2.8 & 7.8 & 50.3   \\
				\hline
				gka3c & -6665 & 60 & 40\% & -6849.9 & -6730.7 & -6703.1 & \textbf{-6688.1} & (\textbf{87.5\%}, 0.3\%) & 6.1 & 9.3 & 15.9 & 62.1   \\
				\hline
				gka4c & -7398 & 70 & 30\% & -7647.1 & -7527.7 & -7494.9 & \textbf{-7462.8} & (\textbf{74.0\%}, 0.9\%) & 13.1 & 18.4 & 24.6 & 88.1   \\
				\hline
				gka5c & -7362 & 80 & 20\% & -7684.5 & -7543.7 & -7474.6 & \textbf{-7412.8} & (\textbf{84.2\%}, 0.7\%) & 15.1 & 27.7 & 40.3 & 112.8   \\
				\hline
				gka6c & -5824 & 90 & 10\% & -6065.8 & -5932.2 & -5869.7 & \textbf{-5847.4} & (\textbf{90.3\%}, 0.4\%) & 10.0 & 11.0 & 19.0 & 57.4   \\
				\hline
				gka7c & -7225 & 100 & 10\% & -7422.7 & -7297.8 & -7264.3 & \textbf{-7248.7} & (\textbf{88.0\%}, 0.3\%) & 12.4 & 13.9 & 22.1 & 55.6   \\
				\hline
				gka1d & -6333 & 100 & 10\% & -6592.7 & -6475.3 & -6403.1 & \textbf{-6369.6} & (\textbf{85.9\%}, 0.6\%) & 11.4 & 13.4 & 29.1 & 71.3   \\
				\hline
				gka2d & -6579 & 100 & 20\% & -7234.2 & -6980.5 & -6897.9 & \textbf{-6811.6} & (\textbf{64.5\%}, 3.5\%) & 42.3 & 70.8 & 70.6 & 193.7   \\
				\hline
				gka3d & -9261 & 100 & 30\% & -9963.0 & -9686.2 & -9591.7 & \textbf{-9523.6} & (\textbf{62.6\%}, 2.8\%) & 164.8 & 200.4 & 262.7 & 330.0   \\
				\hline
				gka4d & -10727 & 100 & 40\% & -11592.5 & -11303.3 & -11175.4 & \textbf{-11096.5} & (\textbf{57.3\%}, 3.4\%) & 302.2 & 259.1 & 191.8 & 387.7   \\
				\hline
				gka5d & -11626 & 100 & 50\% & -12632.1 & -12381.6 & -12274.7 & -12185.0 & (\textbf{44.4\%}, 4.8\%) & 324.3 & 256.3 & 294.3 & 380.2   \\
				\hline
				gka6d & -14207 & 100 & 60\% & -15235.3 & -14938.2 & -14834.9 & \textbf{-14720.2} & (\textbf{50.1\%}, 3.6\%) & 236.6 & 239.7 & 221.9 & 437.9   \\
				\hline
				gka7d & -14476 & 100 & 70\% & -15672.0 & -15413.2 & -15267.6 & -15173.6 & (\textbf{41.7\%}, 4.8\%) & 138.8 & 225.9 & 150.0 & 314.6   \\
				\hline
				gka8d & -16352 & 100 & 80\% & -17353.3 & -17011.5 & -16887.6 & \textbf{-16794.3} & (\textbf{55.8\%}, 2.7\%) & 271.5 & 277.9 & 291.6 & 408.6   \\
				\hline
				gka9d & -15656 & 100 & 90\% & -17010.9 & -16652.0 & -16513.3 & -16409.6 & (\textbf{44.4\%}, 4.8\%) & 390.5 & 419.8 & 367.0 & 513.5   \\
				\hline
				gka10d & -19102 & 100 & 100\% & -20421.4 & -20121.7 & -19974.1 & -19863.8 & (\textbf{44.3\%}, 4.0\%) & 77.8 & 83.4 & 130.2 & 244.8   \\
				\hline
				qplib0032 & 10.1 & 100 & 99\% & -19751 & -16491 & -15962 & -15440 & (21.8\%, 152971.3\%) & 18.1 & 19.4 & 37.0 & 94.7   \\
				\hline
				qplib0067 & -110942 & 80 & 89\% & -116480 & -112923 & -112615 & \textbf{-112478} & (\textbf{72.3\%}, 1.4\%) & 6.2 & 11.1 & 21.9 & 158.3   \\
				\hline
				qplib0633 & 79.6 & 75 & 99\% & 70.9 & 74.0 & 75.1 & \textbf{75.7} & (\textbf{55.2\%}, 4.9\%) & 2.9 & 10.1 & 27.1 & 140.0   \\
				\hline
				qplib2512 & 135028 & 100 & 77\% & -441284 & -125060 & 27898 & \textbf{8290}9 & (\textbf{91.0\%}, 38.6\%) & 18.6 & 19.9 & 53.4 & 278.6   \\
				\hline
				qplib3762 & -296 & 90 & 28\% & -345.6 & -330.8 & -319.9 & \textbf{-309.5} & (\textbf{72.8\%}, 4.6\%) & 6.3 & 18.1 & 50.7 & 183.4   \\
				\hline
				qplib5935 & 4758 & 100 & 99\% & 67494 & 40148 & 36842 & \textbf{28812} & (\textbf{61.7\%}, 505.5\%) & 12.8 & 39.4 & 259.0 & 1745.3   \\
				\hline
				qplib5944 & 1829 & 100 & 99\% & 66934 & 27437 & 23142 & \textbf{19784} & (\textbf{72.4\%}, 981.7\%) & 15.6 & 182.1 & 2304.3 & 13204.6   \\
				\hline
			\end{tabular}
		}
	\end{table}

	\subsubsection*{Max-Cliq instances}
	We take the same graphs as the ones considered in the Max-Cut instances. 
	Some instances share the same adjacency matrix with different weights, in which case we delete these repeated graphs. LB denotes the lower bound of a given instance, computed by $10^6$ random samples. By contrast with the strategy used for the Max-Cut instances, we use sublevel relaxations with level 2 and depth 0, 20, 40, 60, respectively. 
	From Table \ref{maxcliq} we see that the sublevel relaxation yields large improvement compared to Shor's relaxation.  
	The Max-Cliq problem remains hard to solve as emphasized by the large relative gap, ranging from 662.5\% to 3660\%.
	
	\begin{table}[ht]
		\scriptsize
		\caption{Results obtained with sublevel relaxations of Max-Cliq problems.}\label{maxcliq}
		\makebox[1 \textwidth][c]{       
			\begin{tabular}{|c|c|c|c|ccccc|cccc|}
				\hline
				& \multirow{2}{*}{LB} & \multirow{2}{*}{nVar} & \multirow{2}{*}{Density} & \multicolumn{9}{c|}{Sublevel relaxation, $l$ = 2, $q$ = 0/20/40/60 (depth 0 = Shor)} \\
				\cline{5-13}
				& &&& \multicolumn{5}{c|}{upper bounds (RI, RG)} & \multicolumn{4}{c|}{solving time (s)} \\
				\hline
				g05\_60 & 0.8 & 60 & 50\% & 29.9 & 19.3 & 8.3 & \textbf{6.1} & (\textbf{81.8\%}, 662.5\%) & 0.6 & 1.8 & 3.6 & 2.4   \\
				\hline
				g05\_80 & 0.9 & 80 & 50\% & 39.9 & 29.1 & 20.1 & \textbf{8.9} & (\textbf{79.5\%}, 888.9\%) & 2.8 & 7.4 & 7.5 & 8.3   \\
				\hline
				g05\_100 & 0.8 & 100 & 50\% & 50.0 & 39.1 & 28.9 & \textbf{18.4} & (\textbf{64.2\%}, 2200.0\%) & 6.5 & 30.9 & 19.1 & 33.0   \\
				\hline
				pm1d\_80 & 1.0 & 80 & 99\% & 78.2 & 57.5 & 37.6 & \textbf{17.9} & (\textbf{78.1\%}, 1690.0\%) & 2.3 & 5.4 & 7.6 & 4.4   \\
				\hline
				pm1d\_100 & 1.0 & 100 & 99\% & 98.0 & 77.2 & 57.5 & \textbf{37.6} & (\textbf{62.3\%}, 3660\%) & 5.6 & 12.5 & 22.8 & 17.3   \\
				\hline
				pm1s\_80 & 0.7 & 80 & 10\% & 8.9 & 6.2 & 4.6 & \textbf{4.6} & (\textbf{52.1\%}, 557.1\%) & 2.6 & 6.1 & 6.4 & 9.0   \\
				\hline
				pw01\_100 & 0.6 & 100 & 10\% & 10.6 & 8.2 & 5.9 & \textbf{5.4} & (\textbf{51.8\%}, 800.0\%) & 7.5 & 30.7 & 20.0 & 29.6   \\
				\hline
				pw05\_100 & 0.8 & 100 & 50\% & 49.8 & 39.7 & 28.9 & \textbf{18.9} & (\textbf{63.0\%}, 2262.5\%) & 7.6 & 21.9 & 24.0 & 26.5   \\
				\hline
				pw09\_100 & 1.0 & 100 & 90\% & 89.2 & 70.2 & 51.9 & \textbf{34.0} & (\textbf{62.5\%}, 3300\%) & 8.5 & 15.0 & 33.2 & 28.9   \\
				\hline
			\end{tabular}
		}
	\end{table}
	
	\subsubsection*{QCQP instances}
	We take the MIQCP instances from the Biq-Mac library with size larger or equal than 50, then add one dense quadratic constraint $||\mathbf{x}||_2^2 = 1$, and relax the integer bound constraints $\mathbf{x} \in \{0,1\}^n$ to linear bound constraints $\mathbf{x} \in [0,1]^n$. UB denotes the upper bound obtained by selecting the minimum value over $10^6$ random evaluations. 
	
	We also select some instances and their solutions from the QPLIB library with ID 1535, 1661, 1675, 1703 and 1773. These instances have more than one quadratic constraint and involve linear constraints.
	
	Table \ref{qc_summary} is a summary of basic information as well as the number of quadratic, linear, and bound constraints of the instances from the QPLIB library.

	\begin{table}[H]
		\scriptsize
		\caption{Summary of the basic information and constraint structure of the QCQP instances from the QPLIB library.} \label{qc_summary}
		\makebox[1 \textwidth][c]{       
			\begin{tabular}{|c|c|c|c|c|c|c|}
				\hline
				& nVar & Density & nQuad & nLin & nBound \\
				\hline
				qplib1535 & 60 & 94\% & 60 & 6 & 60 \\
				\hline
				qplib1661 & 60 & 95\% & 1 & 12 & 60 \\
				\hline
				qplib1675 & 60 & 49\% & 1 & 12 & 60 \\
				\hline
				qplib1703 & 60 & 98\% & 30 & 6 & 60 \\
				\hline
				qplib1773 & 60 & 95\% & 1 & 6 & 60 \\
				\hline
			\end{tabular}
		}
	\end{table}

	In Table \ref{qc}, we show the lower bounds and running time obtained by the sublevel relaxation with depth 1 for the instances from the QPLIB library, 10 for the instances adapted from the Biq-Mac library, and level 0, 4, 6, 8, respectively. We see that the sublevel relaxation yields a uniform improvement compared to Shor's relaxation. However, for the QCQP problems adapted from the MIQCP instances, it is very hard to find the exact optimal solution as the relative gap varies from 60.5\% to 77.0\%. 
	This is in deep contrast with the instances from the QPLIB library which are relatively easier to solve as the relative gap varies from 9.4\% to 13.8\%.

	\begin{table}[H]
		\scriptsize
		\caption{Results obtained with sublevel relaxations of QCQP problems.} \label{qc}
		\makebox[1 \textwidth][c]{       
			\begin{tabular}{|c|c|c|c|ccccc|cccc|}
				\hline
				& \multirow{2}{*}{Sol./UB} & \multirow{2}{*}{nVar} & \multirow{2}{*}{Density} & \multicolumn{9}{c|}{Sublevel relaxation, $l$ = 0, 4, 6, 8, $q$ = 1, 10 (level 0 = Shor)} \\
				\cline{5-13}
				&&&& \multicolumn{5}{c|}{lower bounds (RI, RG)} & \multicolumn{4}{c|}{solving time (s)} \\
				\hline
				bqp50-1 & -99 & 50 & 10\% & -215.7 & -195.4 & -180.6 & -172.5 & (37.0\%, 74.2\%) & 0.8 & 2.4 & 12.3 & 88.5  \\
				\hline
				bqp100-1 & -67.2 & 100 & 10\% & -323.1 & -304.7 & -296.1 & -290.0 & (12.9\%, 331.5\%) & 21.4 & 22.7 & 56.9 & 249.6   \\
				\hline
				gka1a & -109.5 & 50 & 10\% & -241.8 & -224.1 & -219.8 & -213.8 & (21.2\%, 95.3\%) & 0.8 & 1.9 & 10.0 & 65.6   \\
				\hline
				gka2a & -140.7 & 60 & 10\% & -275.3 & -260.9 & -258.7 & -251.6 & (17.6\%, 78.8\%) & 1.7 & 3.8 & 16.8 & 126.9   \\
				\hline
				gka3a & -143.2 & 70 & 10\% & -300.0 & -284.6 & -278.8 & -275.0 & (15.9\%, 92.0\%) & 3.6 & 9.1 & 23.1 & 121.1   \\
				\hline
				gka4a & -126.2 & 80 & 10\% & -311.0 & -288.3 & -282.9 & -280.0 & (16.8\%, 121.9\%) & 6.8 & 7.9 & 25.4 & 160.0   \\
				\hline
				gka5a & -180.2 & 50 & 20\% & -351.8 & -319.3 & -306.4 & -299.1 & (30.7\%, 66.0\%) & 0.7 & 2.8 & 15.1 & 75.6   \\
				\hline
				gka8a & -122.5 & 100 & 62.5\% & -320.1 & -306.8 & -302.1 & -299.5 & (10.4\%, 144.5\%) & 21.3 & 23.5 & 72.5 & 232.0   \\
				\hline
				gka4b & -63 & 50 & 100\% & -381.4 & -326.2 & -302.4 & -280.8 & (31.6\%, 345.7\%) & 0.7 & 1.8 & 17.4 & 79.2   \\
				\hline
				gka5b & -63 & 60 & 100\% & -446.8 & -377.2 & -348.9 & -327.4 & (31.1\%, 419.7\%) & 1.1 & 4.2 & 19.5 & 117.2   \\
				\hline
				gka6b & -63 & 70 & 100\% & -496.6 & -409.9 & -385.9 & -366.8 & (29.9\%, 482.2\%) & 3.2 & 5.3 & 17.4 & 118.1   \\
				\hline
				gka7b & -63 & 80 & 100\% & -518.3 & -447.1 & -421.6 & -404.3 & (25.0\%, 541.7\%) & 5.9 & 9.5 & 21.6 & 170.3   \\
				\hline
				gka8b & -63 & 90 & 100\% & -534.5 & -472.7 & -449.0 & -430.1 & (22.1\%, 582.7\%) & 10.9 & 16.0 & 36.0 & 148.9   \\
				\hline
				gka9b & -63 & 100 & 100\% & -573.0 & -501.3 & -477.0 & -455.8 & (23.0\%, 623.5\%) & 19.7 & 36.6 & 42.5 & 191.4   \\
				\hline
				gka10b & -63 & 125 & 100\% & -639.4 & -569.7 & -553.7 & -533.6 & (18.4\%, 747.0\%) & 80.1 & 82.1 & 110.9 & 410.1   \\
				\hline
				gka2c & -159.1 & 50 & 60\% & -290.0 & -269.3 & -261.6 & -255.4 & (26.4\%, 60.5\%) & 0.8 & 2.5 & 11.2 & 80.8   \\
				\hline
				gka3c & -126.3 & 60 & 40\% & -271.2 & -240.2 & -235.4 & -231.3 & (27.5\%, 83.1\%) & 1.7 & 4.4 & 16.0 & 103.1   \\
				\hline
				gka4c & -123.0 & 70 & 30\% & -292.7 & -263.7 & -254.4 & -247.9 & (26.4\%, 101.5\%) & 3.0 & 6.5 & 19.7 & 155.6   \\
				\hline
				gka5c & -114.0 & 80 & 20\% & -239.1 & -225.9 & -223.2 & -220.4 & (14.9\%, 93.3\%) & 10.6 & 9.6 & 32.3 & 166.5   \\
				\hline
				gka6c & -100 & 90 & 10\% & -198.8 & -190.8 & -186.7 & -182.4 & (16.6\%, 82.4\%) & 12.3 & 15.9 & 33.1 & 216.5   \\
				\hline
				gka7c & -100 & 100 & 10\% & -225.8 & -213.7 & -210.5 & -208.5 & (13.8\%, 108.5\%) & 21.4 & 28.2 & 63.4 & 323.1   \\
				\hline
				gka1d & -75 & 100 & 10\% & -197.9 & -182.5 & -177.0 & -174.5 & (19.0\%, 132.7\%) & 19.3 & 25.9 & 64.8 & 243.9   \\
				\hline
				gka2d & -87.2 & 100 & 20\% & -259.6 & -242.2 & -233.8 & -229.5 & (17.5\%, 163.2\%) & 23.5 & 28.2 & 61.1 & 254.3   \\
				\hline
				gka3d & -88.1 & 100 & 30\% & -304.0 & -281.6 & -274.1 & -267.5 & (16.9\%, 203.6\%) & 26.2 & 28.9 & 49.2 & 278.4   \\
				\hline
				gka4d & -105.5 & 100 & 40\% & -375.2 & -340.1 & -326.0 & -317.5 & (21.4\%, 201.0\%) & 21.6 & 21.9 & 53.2 & 270.7   \\
				\hline
				gka5d & -131.9 & 100 & 50\% & -383.6 & -351.5 & -341.5 & -332.3 & (20.4\%, 152.0\%) & 20.4 & 22.7 & 41.3 & 257.1   \\
				\hline
				gka6d & -137.7 & 100 & 60\% & -443.1 & -400.0 & -391.0 & -378.9 & (21.0\%, 175.2\%) & 23.8 & 23.3 & 48.6 & 254.2   \\
				\hline
				gka7d & -156.3 & 100 & 70\% & -453.9 & -421.4 & -406.6 & -397.4 & (19.0\%, 154.3\%) & 21.4 & 22.5 & 71.7 & 217.8   \\
				\hline
				gka8d & -147.6 & 100 & 80\% & -488.0 & -441.1 & -423.3 & -414.2 & (21.7\%, 180.6\%) & 21.7 & 25.0 & 47.0 & 232.8   \\
				\hline
				gka9d & -179.6 & 100 & 90\% & -539.7 & -487.7 & -469.2 & -456.8 & (23.0\%, 154.3\%) & 20.6 & 21.5 & 45.1 & 222.2   \\
				\hline
				gka10d & -187.0 & 100 & 100\% & -552.4 & -505.7 & -491.8 & -478.4 & (20.3\%, 155.8\%) & 23.2 & 24.6 & 56.9 & 196.5   \\
				\hline
				qplib1535 & -11.6 & 60 & 94\% & -13.9 & -13.5 & -13.3 & -13.2 & (30.4\%, 13.8\%) & 1.4 & 4.3 & 13.7 & 99.2   \\
				\hline
				qplib1661 & -16.0 & 60 & 95\% & -18.4 & -18.1 & -17.8 & -17.5 & (37.5\%, 9.4\%) & 1.4 & 3.0 & 14.7 & 96.4   \\
				\hline
				qplib1675 & -75.7 & 60 & 49\% & -93.1 & -87.0 & -85.2 & \textbf{-83.8} & (\textbf{53.4\%}, 10.7\%) & 1.0 & 4.2 & 19.2 & 147.8   \\
				\hline
				qplib1703 & -132.8 & 60 & 98\% & -152.8 & -147.0 & -145.2 & \textbf{-143.5} & (\textbf{46.5\%}, 8.06\%) & 1.2 & 4.2 & 20.8 & 109.3   \\
				\hline
				qplib1773 & -14.6 & 60 & 95\% & -17.3 & -16.8 & -16.6 & -16.4 & (33.3\%, 12.3\%) & 1.1 & 4.0 & -14.0 & 89.3  \\
				\hline
			\end{tabular}
		}
	\end{table}

	\subsection{Examples from deep learning}
	\subsubsection*{Lipschitz constant estimation}
	We generate random 1-hidden layer neural networks with parameters $\mathbf{A}, \mathbf{b},\mathbf{c}$. We denote by net\_1\_$n$ the instances of 1-hidden layer networks of size $n$, and compute the upper bounds corresponding to $\epsilon = 0.1, 10$, by the sublevel relaxations of depth 1 and level 0, 4, 6, 8, respectively. 
	
	We see in Table \ref{globallip} that we have a relatively high improvement ratio and low gap for the global case, while Table \ref{locallip} ,dedicated to the local case, shows that the improvement ratio is decreasing and the gap is increasing. 
	The underlying rationale is that local Lipschitz constants of neural networks are harder to estimate than the global ones.
	
	\begin{table}[H]
		\scriptsize
		\caption{Results obtained with sublevel relaxation of Lipschitz constant problems, $\epsilon = 10$.} \label{globallip}
		\makebox[1 \textwidth][c]{       
			\begin{tabular}{|c|c|c|ccccc|cccc|}
				\hline
				& \multirow{2}{*}{Sol./LB} & \multirow{2}{*}{nVar} & \multicolumn{9}{c|}{Sublevel relaxation, lv = 0/4/6/8, p = 1 (level 0 = Shor)} \\
				\cline{4-12}
				& && \multicolumn{5}{c|}{upper bounds (RI, RG)} & \multicolumn{4}{c|}{solving time (s)} \\
				\hline
				net\_1\_5 & 0.38 & 15 & 0.44 & 0.39 & \textbf{0.38} & 0.38 & (\textbf{100\%}, 0\%) & 0.02 & 0.61 & 1.77 & 7.01 \\
				\hline
				net\_1\_10 & 0.69 & 30 & 0.72 & 0.70 & \textbf{0.69} & 0.69 & (\textbf{100\%}, 0\%) & 0.10 & 0.90 & 4.18 & 26.05 \\
				\hline
				net\_1\_15 & 1.72 & 45 & 1.86 & 1.81 & 1.76 & \textbf{1.73} & (\textbf{92.86\%}, 0.58\%) & 0.35 & 2.19 & 10.55 & 69.66 \\
				\hline
				net\_1\_20 & 2.68 & 60 & 2.88 & 2.83 & 2.77 & \textbf{2.75} & (\textbf{65.00\%}, 2.61\%) & 1.17 & 4.24 & 15.56 & 60.59 \\
				\hline
				net\_1\_25 & 3.56 & 75 & 3.83 & 3.74 & 3.69 & \textbf{3.68} & (\textbf{55.56\%}, 3.37\%) & 2.79 & 8.72 & 29.38 & 166.35 \\
				\hline
				net\_1\_30 & 5.60 & 90 & 6.16 & 6.11 & 6.08 & 6.06 & (17.86\%, 8.21\%) & 8.45 & 11.68 & 33.29 & 220.39 \\
				\hline
				net\_1\_35 & 7.77 & 105 & 8.92 & 8.79 & 8.73 & 8.66 & (22.61\%, 11.455\%) & 16.69 & 26.55 & 74.28 & 267.19 \\
				\hline
				net\_1\_40 & 7.40 & 120 & 9.07 & 8.97 & 8.86 & 8.78 & (17.37\%, 18.65\%) & 33.19 & 56.15 & 116.37 & 333.65 \\
				\hline
			\end{tabular}
		}
	\end{table}

	\begin{table}[H]
		\scriptsize
		\caption{Results obtained with sublevel relaxations of Lipschitz constant problems, $\epsilon = 0.1$.} \label{locallip}
		\makebox[1 \textwidth][c]{       
			\begin{tabular}{|c|c|c|ccccc|cccc|}
				\hline
				& \multirow{2}{*}{Sol./LB} & \multirow{2}{*}{nVar} & \multicolumn{9}{c|}{Sublevel relaxation, lv = 0/4/6/8, p = 1 (level 0 = Shor)} \\
				\cline{4-12}
				& && \multicolumn{5}{c|}{upper bounds (RI, RG)} & \multicolumn{4}{c|}{solving time (s)} \\
				\hline
				net\_1\_5 & 0.247 & 15 & 0.251 & 0.251 & \textbf{0.247} & 0.247 & (\textbf{100\%}, 0\%) & 0.04 & 0.34 & 1.25 & 6.72 \\
				\hline
				net\_1\_10 & 0.581 & 30 & 0.610 & 0.608 & 0.606 & 0.605 & (17.2\%, 4.13\%) & 0.18 & 0.84 & 4.65 & 38.56 \\
				\hline
				net\_1\_15 & 1.384 & 45 & 1.449 & 1.441 & 1.441 & 1.435 & (21.54\%, 3.68\%) & 0.42 & 1.43 & 7.29 & 60.27 \\
				\hline
				net\_1\_20 & 1.73 & 60 & 2.23 & 2.22 & 2.20 & 2.19 & (8.00\%, 26.59\%) & 4.21 & 3.82 & 13.11 & 83.14 \\
				\hline
				net\_1\_25 & 2.03 & 75 & 2.73 & 2.67 & 2.65 & 2.64 & (12.86\%, 30.05\%) & 4.79 & 7.08 & 23.31 & 134.50 \\
				\hline
				net\_1\_30 & 4.10 & 90 & 5.09 & 5.07 & 5.06 & 5.04 & (5.05\%, 22.93\%) & 19.10 & 13.60 & 28.29 & 146.17 \\
				\hline
				net\_1\_35 & 5.84 & 105 & 7.12 & 7.08 & 7.07 & 7.03 & (7.03\%, 20.38\%) & 56.46 & 28.95 & 47.06 & 192.31 \\
				\hline
				net\_1\_40 & 5.02 & 120 & 7.30 & 7.21 & 7.15 & 7.07 & (10.09\%, 40.84\%) & 144.28 & 58.01 & 80.27 & 254.22 \\
				\hline
			\end{tabular}
		}
	\end{table}

	\subsubsection*{Certification instances}
	We use the same network net\_1\_$n$ as the one generated for the above Lipschitz problems, and compute the upper bounds corresponding to $\epsilon = 0.1, 10$, by the sublevel relaxations of depth 1 and level 0, 4, 6, 8, respectively. 
	
	As for the Lipschitz problem, Table \ref{globalcert} and \ref{localcert} indicate that in the local case it is much harder to improve and find the exact optimal solution than in the global case. Furthermore, the difficulty of the problem also increases with the dimension. 
	When the number of variables gets larger, the improvement ratio decreases while the relative gap increases.
	
	\begin{table}[H]
		\scriptsize
		\caption{Results obtained with sublevel relaxations of certification problems, $\epsilon = 10$.} \label{globalcert}
		\makebox[1 \textwidth][c]{       
			\begin{tabular}{|c|c|c|ccccc|cccc|}
				\hline
				& \multirow{2}{*}{Sol./LB} & \multirow{2}{*}{nVar} & \multicolumn{9}{c|}{Sublevel relaxation, lv = 0/4/6/8, p = 1 (level 0 = Shor)} \\
				\cline{4-12}
				& && \multicolumn{5}{c|}{upper bounds (RI, RG)} & \multicolumn{4}{c|}{solving time (s)} \\
				\hline
				net\_1\_5 & 2.63 & 10 & 3.51 & 3.00 & \textbf{2.74} & 2.74 & (\textbf{87.50\%}, 4.18\%) & 0.01 & 0.26 & 1.79 & 1.83 \\
				\hline
				net\_1\_10 & 3.49 & 20 & 4.88 & 4.69 & 4.60 & 4.48 & (28.78\%, 28.37\%) & 0.06 & 0.99 & 4.23 & 30.36 \\
				\hline
				net\_1\_15 & 5.61 & 30 & 8.20 & 8.10 & 7.84 & 7.41 & (30.50\%, 32.09\%) & 0.19 & 1.02 & 6.93 & 40.31 \\
				\hline
				net\_1\_20 & 9.24 & 40 & 16.48 & 16.03 & 15.75 & 15.48 & (13.81\%, 67.53\%) & 0.60 & 2.35 & 9.69 & 67.31 \\
				\hline
				net\_1\_25 & 14.40 & 50 & 26.68 & 26.28 & 25.89 & 25.57 & (9.04\%, 77.57\%) & 2.24 & 5.67 & 17.40 & 66.99 \\
				\hline
				net\_1\_30 & 17.22 & 60 & 38.06 & 37.72 & 36.82 & 35.89 & (10.41\%, 108.42\%) & 5.08 & 14.32 & 24.02 & 102.93 \\
				\hline
				net\_1\_35 & 26.71 & 70 & 59.18 & 58.64 & 57.78 & 57.39 & (5.51\%, 114.86\%) & 10.69 & 25.79 & 40.96 & 136.35 \\
				\hline
				net\_1\_40 & 22.94 & 80 & 57.59 & 56.08 & 54.69 & 54.18 & (9.84\%, 136.18\%) & 23.22 & 44.74 & 65.04 & 146.91 \\
				\hline
				net\_1\_45 & 22.57 & 90 & 57.56 & 56.34 & 55.57 & 54.68 & (8.23\%, 142.27\%) & 44.67 & 85.38 & 107.12 & 186.94 \\
				\hline
				net\_1\_50 & 27.34 & 100 & 73.59 & 72.10 & 71.19 & 69.92 & (7.94\%, 155.74\%) & 81.61 & 144.61 & 165.18 & 333.42 \\
				\hline
			\end{tabular}
		}
	\end{table}
	
	\begin{table}[H]
		\scriptsize
		\caption{Results obtained with sublevel relaxations of certification problems, $\epsilon = 0.1$.} \label{localcert}
		\makebox[1 \textwidth][c]{       
			\begin{tabular}{|c|c|c|ccccc|cccc|}
				\hline
				& \multirow{2}{*}{Sol./LB} & \multirow{2}{*}{nVar} & \multicolumn{9}{c|}{Sublevel relaxation, lv = 0/4/6/8, p = 1 (level 0 = Shor)} \\
				\cline{4-12}
				& && \multicolumn{5}{c|}{upper bounds (RI, RG)} & \multicolumn{4}{c|}{solving time (s)} \\
				\hline
				net\_1\_5 & 0.190 & 10 & 0.191 & 0.191 & 0.191 & 0.191 & (0.00\%, 0.53\%) & 0.01 & 0.13 & 0.74 & 0.73 \\
				\hline
				net\_1\_10 & 0.021 & 20 & 0.025 & 0.025 & 0.025 & 0.024 & (25.00\%, 14.29\%) & 0.15 & 0.51 & 3.16 & 17.68 \\
				\hline
				net\_1\_15 & 0.027 & 30 & 0.053 & 0.053 & 0.053 & 0.053 & (0.00\%, 96.30\%) & 0.17 & 0.55 & 2.82 & 17.01 \\
				\hline
				net\_1\_20 & 0.269 & 40 & 0.299 & 0.299 & 0.299 & 0.298 & (3.33\%, 10.78\%) & 0.79 & 2.36 & 7.83 & 46.08 \\
				\hline
				net\_1\_25 & -0.104 & 50 & -0.025 & -0.028 & -0.031 & -0.031 & (7.59\%, 70.19\%) & 2.37 & 5.77 & 13.22 & 54.81 \\
				\hline
				net\_1\_30 & 0.669 & 60 & 0.810 & 0.807 & 0.806 & 0.803 & (4.96\%, 20.03\%) & 6.34 & 10.06 & 21.85 & 84.82 \\
				\hline
				net\_1\_35 & 0.825 & 70 & 1.107 & 1.107 & 1.107 & 1.107 & (0.00\%, 34.18\%) & 12.61 & 18.44 & 34.07 & 102.26 \\
				\hline
				net\_1\_40 & 0.741 & 80 & 0.949 & 0.943 & 0.942 & 0.940 & (4.33\%, 26.86\%) & 34.70 & 52.44 & 56.64 & 163.67 \\
				\hline
				net\_1\_45 & 0.265 & 90 & 0.603 & 0.602 & 0.600 & 0.599 & (1.18\%, 126.04\%) & 55.69 & 89.39 & 115.78 & 200.94 \\
				\hline
				net\_1\_50 & 0.614 & 100 & 0.920 & 0.919 & 0.916 & 0.914 & (1.96\%, 48.86\%) & 105.68 & 177.41 & 179.53 & 205.58 \\
				\hline
			\end{tabular}
		}
	\end{table}
	
	The two latter examples show us that finding the guaranteed bounds for optimization problems arising from deep learning is much harder than the usual sparse problems coming from the classical optimization literature. 
	Hence, it remains a big challenge to adapt our approach to large real networks, involving a large number of variables and more complicated structures such as  convolutional or max-pooling layers.
	
	\section{Conclusion}
	In this paper, we propose a new semidefinite programming hierarchy based on the standard dense and sparse Lasserre's hierarchies. 
	This hierarchy provides a wider choice of intermediate relaxation levels, lying between the $d$-th and $(d+1)$-th order relaxations in  Lasserre's hierarchy. With this technique, we are able to solve problems where the standard relaxations are untractable.  
	Our experimental results demonstrate that the sublevel relaxation often allows one to compute more accurate bounds by comparison with existing frameworks such as Shor's relaxation or term sparsity, in particular for dense problems.
	
	Sublevel relaxations offer a large choice of parameters tuning, as one can select   the level, depth, and subsets for each relaxation.
	We can benefit from this to potentially perform better that state-of-the-art methods.
	However, the flexibility of our approach also comes together with a drawback since the more flexible it is, the more difficult for the users it is to tune the parameters. 
	One important and interesting future topic would be to design an algorithm that searches for the optimal level, depth and subsets in sublevel relaxations.
	
	\section*{Acknowledgement}
	This work has benefited from the Tremplin ERC Stg Grant ANR-18-ERC2-0004-01 (T-COPS project), the European Union's Horizon 2020 research and innovation programme under the Marie Sklodowska-Curie Actions, grant agreement 813211 (POEMA) as well as from the AI Interdisciplinary Institute ANITI funding, through the French ``Investing for the Future PIA3'' program under the Grant agreement n$^{\circ}$ANR-19-PI3A-0004. The third author was supported by the FMJH Program PGMO (EPICS project) and  EDF, Thales, Orange et Criteo.
	The fourth author acknowledge the support of Air Force Office of Scientific Research, Air Force Material Command, USAF, under grant numbers FA9550-19-1-7026, FA9550-18-1-0226, and ANR MasDol.

	\bibliographystyle{plain}
	\bibliography{sublevel}
	
\end{document}